\documentclass{article}

\usepackage{cmap}
\usepackage[T2A]{fontenc}
\usepackage[cp866]{inputenc}
\usepackage[english,russian]{babel}
\usepackage[tbtags]{amsmath}
\usepackage{amsfonts,amssymb,mathrsfs,amscd}

\usepackage{hypbmsec} 

\newtheorem{probl}{Задача}   
\newtheorem{probll}{Задача}   
\newtheorem{theorem}{Теорема}
  
\newtheorem{prop}{Предложение}

\newtheorem{proof}{Доказательство}   

\renewcommand{\theenumi}{\arabic{enumi})}
\renewcommand{\labelenumi}{\theenumi}

\def\CC{\mathbb C}
\def\RR{\mathbb R}

\def\supp{\operatorname{supp}}
\let\myo\overline\let\myh\widehat
\let\le\leqslant\let\ge\geqslant

\let\le=\leqslant
\let\ge=\geqslant
\def\bad{\spaceskip=0.33emplus0.6emminus0.15em\immediate\write5{\string\bad}}
\begin{document}

\author{С.\,П.~Суетин}
\date{16.10.2013}

\title{Формулы следов для некоторого класса операторов Якоби}

\maketitle

\markright{Формулы следов для некоторого класса операторов
Якоби}

\footnotetext[0]{Работа выполнена при поддержке
Российского фонда фундаментальных исследований (грант
\No~11-01-00330а) и Программы
поддержки ведущих научных школ РФ (грант
\No~НШ-4664.2012.1).}

\begin{abstract}
В работе изучается класс операторов Якоби таких, что каждый оператор 
порождается единичной борелевской мерой с носителем, состоящим из конечного 
числа отрезков на вещественной прямой $\mathbb R$ и конечного числа точек в 
$\mathbb C$, расположенных вне выпуклой оболочки этих отрезков и симметрично 
относительно $\mathbb R$. В таком классе операторов получены асимптотика
диагональной функции Грина и формулы следов для последовательностей 
коэффициентов, соответствующих заданному оператору.

Библиография: 34 названия.
\end{abstract}
\section{Введение}
\label{s1}

\subsection{}
\label{s1.1}
Пусть функция $f=\myh{\mu}$, где
\begin{equation}
\label{1}
\myh{\mu}(z)=\int\frac{d\mu(t)}{z-t},\quad z\in\myo\CC\setminus\supp{\mu},
\end{equation}
-- преобразование Коши положительной меры $\mu$  с компактынм носителем в
$\CC$, состоящим из конечного числа отрезков, расположенных на вещественной
прямой $\RR$, и конечного числа точек вне этих отрезков, расположенных
симметрично относительно $\RR$. А именно,
$\supp\mu=E\sqcup\{\zeta_1,\dots,\zeta_m\}$, где $E=\bigsqcup_{j=1}^{g+1}E_j$,
$E_j=[e_{2j-1},e_{2j}]$ и множество $\{\zeta_1,\dots,\zeta_m\}$ расположено
симметрично относительно вещественной прямой и не пересекается с $\myh{E}$
-- выпуклой оболочкой $E$. В дальнейшем предполагается, что на $E$ мера
$\mu$ абсолютно неперерывна относительно меры Лебега и имеет следующий вид
\begin{equation}
\label{2}
\frac{d\mu}{dx}=\frac1\pi\frac{\rho(x)}{\sqrt{-h(x)}}>0,\quad x\in E,
\end{equation}
где функция
$h(z)=\prod_{j=1}^{2g+2}(z-e_j)$, причем в~области
$D=\widehat{\mathbb C}\setminus E$ выбрана та ветвь
квадратного корня, которая положительна при положительных
значениях аргумента. Тем самым, $\sqrt{-h(x)}>0$ при
$x\in[e_{2g+1},e_{2g+2}]$; на остальных отрезках значение
корня определяется аналитическим продолжением. Весовая
функция $\rho$ предполагается \textit{голоморфной и отличной
от нуля на~$E$}. Вне~$E$ мера~$\mu$ имеет вид
\begin{equation}
\label{3}
d\mu(\zeta)=\sum_{k=1}^m\lambda_k\delta(\zeta-\zeta_k)\,d\zeta,
\end{equation}
где все $\lambda_k>0$, точки
$\zeta_k\in\mathbb C\setminus\widehat E$.
$\widehat E$~-- выпуклая оболочка~$E$.
Тем самым,
\begin{equation}
\label{4}
\widehat\mu(z)=\frac1\pi\int_E\frac{\rho(x)}{z-x}\,\frac{dx}{\sqrt{-h(x)}}
+\sum_{k=1}^m\frac{\lambda_k}{z-\zeta_k}\,,
\qquad z\notin\operatorname{supp}\mu,
\end{equation}
$\operatorname{supp}\mu=S\cup\{\zeta_1,\dots,\zeta_m\}$.

Хорошо известно (см., например,\cite{3}), что при отсутствии
вырождения функция $f$ с помощью функционального аналога классического
алгоритма Евклида разлагается в непрерывную $J$-дробь (или, иначе говоря,
чебышёвскую непрерывную дробь; см.\ \cite{5},~\cite{2}):
\begin{align}
f(z):&=\frac1{z-b_1-f_1(z)}
=\cfrac1{z-b_1-\cfrac{a_1^2}{z-b_2-f_2(z)}}
\label{eq2}
\\
&\simeq\cfrac[l]{1}
{z-b_1-\cfrac[l]{a_1^2}
{z-b_2-\cfrac[l]{a_2^2\qquad}{\ddots} }}\;.
\label{eq3}
\end{align}
Коэффициенты $a_n$ и $b_n$ непрерывной дроби~\eqref{eq3}
строятся непосредственно по коэф\-фициентам $s_k$ разложения
функции $\widehat\mu$ в~ряд Лорана в~бесконечно удаленной
точке $z=\infty$
\begin{equation}
\label{eq4}
\widehat\mu(z)=\sum_{k=0}^\infty\frac{s_k}{z^{k+1}}\,,
\end{equation}
где
$$
s_k=\int \zeta^k\,d\mu(x), \qquad k=0,1,2,\dots,
$$
-- моменты меры~$\mu$. Следовательно, все параметры $a_n,b_n$ -- вещественные.

Классическая теорема Маркова (см.\ 
\cite{7} (3.5, теорема~3.5.4), \cite{8},~\cite{2})
утверждает, что при $\supp{\mu}\subset\RR$ $J$-дробь~\eqref{eq3}
сходится к~функции $\widehat\mu$ равномерно внутри (т.е.\
на компактных подмножествах) области
$\widehat{\mathbb C}\setminus[A,B]$, где $[A,B]$~-- выпуклая
оболочка $\operatorname{supp}\mu$. По функции~$\widehat\mu$
мера может быть восстановлена с~помощью
формулы Стилтьеса--Перрона:
\begin{equation}
\label{eq5}
\frac{\mu(\zeta_2+0)+\mu(\zeta_2-0)}2-\frac{\mu(\zeta_1+0)+\mu(\zeta_1-0)}2
=-\lim_{\varepsilon\to+0}
\frac1{\pi}\int_{\zeta_1}^{\zeta_2}\operatorname{Im}
\widehat\mu(x+i\varepsilon)\,dx
\end{equation}
(здесь интеграл в~представлении~\eqref{eq5} для
$\widehat\mu$ понимается в~смысле Римана--Стил\-тьеса, т.е.\
$\mu$ рассматривается как функция ограниченной вариации).

Однако непосредственно ни алгоритмическая
процедура~\eqref{eq2}, ни формула~\eqref{eq5} не дают
ответа на вопрос о~том, как те или иные свойства
``регулярности'' меры~$\mu$ связаны с~``регулярностью''
последовательностей $a$ и~$b$, точнее~-- с~асимптотическим
поведением коэффициентов $a_n$ и~$b_n$ чебышёвской
непрерывной дроби~\eqref{eq3} при $n\to\infty$.

Наиболее изученным в~этом отношении является
\textit{классический} случай, когда последовательности $a$
и~$b$ таковы, что $a_n$ и $b_n$ имеют пределы. Точнее, пусть
\begin{equation}
\label{eq6}
a_n\to\frac12\,, \qquad b_n\to0.
\end{equation}
Как хорошо известно (см., например, \cite{9},~\cite{3}),
при этом условии
$\operatorname{ess\,supp}\mu=[-1,1]$ и на
$\mathbb R\setminus[-1,1]$ мера~$\mu$ имеет не более чем
счетное число точечных масс, которые могут накапливаться
лишь к~точкам~$\pm1$~-- концам отрезка $\Delta:=[-1,1]$.
В~свою очередь, соотношение~\eqref{eq6}
справедливо \cite{10},~\cite{11},~\cite{8} для мер~$\mu$,
удовлетворяющих условию $\mu'(x)=d\mu/dx>0$ п.в.\ на
$\operatorname{ess\,supp}\mu=\Delta$ и имеющих конечное
число точечных масс на~$\mathbb R$ вне~$\Delta$.

Соотношения~\eqref{eq6} означают, что коэффициенты $a_n$
и~$b_n$ \textit{предельно постоянны}. Эти соотношения,
очевидно, не выполняются, если
в~$E_\mu:=\operatorname{ess\,supp}\mu$ имеется хотя бы одна
лакуна. Случай \textit{предельно периодических} с~периодом
$N\ge2$ коэффициентов $a_n$ и~$b_n$ соответствует
ситуации, когда в~$E_\mu$ имеется $N-1$
\textit{лакуна}\footnote{Предполагается, что все лакуны ``открытые''.}
$(\alpha_j,\beta_j)$, $j=1,2,\dots,N-1$, тем самым
$E_\mu=[-1,1]\setminus\bigcup_{j=1}^{N-1}(\alpha_j,\beta_j)$.
При этом $E_\mu=\{x:T_{2N}(x)\in[-1,1]\}$, где $T_{2N}$~--
полином с~вещественными коэффициентами степени~$2N$. Тем
самым, компакт~$E_\mu$ состоит из $N$~компонент равной
гармонической меры $\omega_j(\infty)=1/N$, $j=1,\dots,N$.
Условия регулярности асимптотического поведения $a_n$
и~$b_n$, вполне аналогичные приведенным выше
условиям для классического
случая, позволяют сделать определенные заключения
о~регулярности меры~$\mu$ и в~предельно периодическом случае
(подробнее см.~\cite{15}).

В~настоящей работе мы обобщаем результаты работы~\cite{Sue07} и
заранее предполагаем, что существенная часть~$E_\mu$
носителя меры~$\mu$ состоит ровно из $N$ непересекающихся
отрезков, а~сама мера на~$E_\mu$ удовлетворяет
определенному условию регулярности (подробнее см.~\S\,\ref{s2}).
Однако при этом мы не налагаем никаких
априорных ограничений на взаимное расположение компонент
компакта~$E_\mu$, в~частности, они могут находиться
``в~общем положении'' (такому случаю соответствуют рационально
независимые гармонические меры
$\omega_1(\infty),\dots,\omega_N(\infty)$). При этом
условии устанавливаются асимптотические формулы для
коэффициентов $a_n$ и~$b_n$ чебышёвской непрерывной
дроби~\eqref{eq3}, ``заменяющие'' классические
соотношения~\eqref{eq6}. Отметим, что из этих формул
вытекает, что $a_n$ и~$b_n$ являются ``предельно
квазипериодичными'' при $n\to\infty$ с~группой периодов
$\omega_1(\infty),\dots,\omega_{N-1}(\infty)$
(см.~\S\,\ref{s2}, теорема~\ref{t2}). Эти предельные соотношения
являются следствием асимптотической формулы для
диагональной функции Грина самосопряженного оператора
Якоби $J=J_{a,b}$, порожденного парой $a$,~$b$, или, эквивалентно,
мерой~$\mu$ (теорема~\ref{t1}).

\subsection{}
\label{s1.2}
Пусть последовательности $a$,~$b$ возникают из~\eqref{2}
$\{\mathbf e_j\}_{j=1}^\infty$~--
стандартный базис в~$\ell^2=\ell^2(\mathbb N)$,
$\langle\,\cdot\,{,}\,\cdot\,\rangle$~-- скалярное
произведение. На этом базисе стандартным образом определим
симметрический оператор Якоби $J\colon\ell^2\,{\to}\,\nobreak\ell^2$,
$J=J_{a,b}=J_\mu$:
\begin{equation}
\label{eq10}
\begin{aligned}
J\mathbf e_1&=b_1\mathbf e_1+a_1\mathbf e_2,
\\
J\mathbf e_n&=a_{n-1}\mathbf e_{n-1}+b_n\mathbf e_n+a_n\mathbf e_{n+1},
\qquad n=2,3,\dotsc.
\end{aligned}
\end{equation}
Хорошо известно~\cite{2}, что при рассматриваемыфх здесь условиях
оператор~$J$ единственным образом
продолжается как самосопряженный оператор на все~$\ell^2$,
при этом $\mu$~-- спектральная мера~$J$, а~спектр
$\sigma(J)$ оператора~$J$ совпадает с~носителем меры~$\mu$:
$\sigma(J)=\operatorname{supp}\mu$. Функция
$m(z):=\langle\mathbf e_1,(J-z)^{-1}\mathbf e_1\rangle=-\widehat\mu(z)$
называется функцией Вейля оператора~$J$,
$G_{n,m}(z):=\langle\mathbf e_n,(J-z)^{-1}\mathbf e_m\rangle$,
$n,m\in\mathbb N$,~-- его функция Грина,
$G_n(z):=G_{n,n}(z)=\langle\mathbf e_n,(J-z)^{-1}\mathbf e_n\rangle$~--
диагональная функция Грина (см.,
например, \cite{14}, \cite{19},~\cite{21}).

Пусть $q_n(x)=k_n\zeta^n+\dotsb$, $k_n>0$, $n=0,1,\dots$, --
последовательность многочленов, определенная рекуррентными
соотношениями
\begin{equation}
\label{eq11}
a_{n-1}q_{n-2}(x)+b_nq_{n-1}(x)+a_nq_n(x)=xq_{n-1}(x),
\qquad n=1,2\dots,
\end{equation}
где $a_0=1$ (это связано с~условием $\mu(\mathbb R)=1$),
$q_{-1}(x)\equiv0$, $q_0(x)\equiv1$. Нетрудно увидеть, что
$\mathbf e_n=q_{n-1}(J)\mathbf e_1$ и в~соответствии со
спектральной теоремой
\begin{align*}
\delta_{kj}
&=\langle\mathbf e_k,\mathbf e_j\rangle
=\bigl\langle q_{k-1}(J)\mathbf e_1,q_{j-1}(J)\mathbf e_1\bigr\rangle
=\bigl\langle q_{j-1}(J)q_{k-1}(J)\mathbf e_1,\mathbf e_1\bigr\rangle
\\
&=\int q_{k-1}(x)q_{j-1}(x)\,d\mu(x).
\end{align*}
Следовательно, полиномы $\{q_n\}_{n=0}^\infty$
ортонормированы относительно спектральной меры~$\mu$. Для
соответствующих многочленов второго рода $p_n(z)$ и
функций второго рода $r_n(z)$ имеем:
\begin{align}
p_n(z)&:=\int\frac{q_n(z)-q_n(x)}{z-x}\,d\mu(x),
\label{eq12}
\\
r_n(z)&:=\int\frac{q_n(x)\,d\mu(x)}{z-x}
=\frac1{q_n(z)}\int\frac{q^2_n(x)\,d\mu(x)}{z-x}\,.
\label{eq13}
\end{align}
Тем самым,
\begin{equation}
\label{eq14}
r_n(z)=(q_n\widehat\mu-p_n)(z)=\frac1{k_nz^{n+1}}+\dotsb,
\qquad z\to\infty.
\end{equation}
Из~\eqref{eq12} и~\eqref{eq13} вытекает, что
последовательности функций $p_n(z)$ и $r_n(z)$ также
удовлетворяют рекуррентным соотношениям~\eqref{eq11}, но
с~другими начальными условиями:
$$
p_{-1}(z)\equiv-1,\quad
p_0(z)\equiv0,\quad r_{-1}(z)\equiv1,\quad
r_0(z)=\widehat\mu(z).
$$

Отметим, что из~\eqref{eq14} вытекает соотношение
\begin{equation}
\widehat\mu(z)-\frac{p_n}{q_n}(z)
=O\biggl(\frac1{z^{2n+1}}\biggr),
\qquad z\to\infty.
\label{eq15}
\end{equation}
Тем самым, рациональная функция~$p_n/q_n$ есть $n$-я
\textit{подходящая дробь} к~чебышёвской непрерывной
дроби~\eqref{eq3} или, иначе говоря, $n$-я
\textit{диагональная аппроксимация Паде} ряда~\eqref{eq4}.

\goodbreak

С~помощью соотношения $\mathbf e_n=q_{n-1}(J)\mathbf e_1$
и спектральной теоремы получаем, что
$$
G_{n,m}(z)=\int\frac{q_{n-1}(x)q_{m-1}(x)}{z-x}\,d\mu(x)
=\begin{cases}
q_{n-1}(z)r_{m-1}(z),&n\le m,
\\
q_{m-1}(z)r_{n-1}(z),&n\ge m.
\end{cases}
$$
Тем самым,
\begin{equation}
\label{eq16}
G_{n+1}(z)=\int\frac{q_n^2(x)}{z-x}\,d\mu(x)=q_n(z)r_n(z),
\qquad z\notin\operatorname{supp}\mu.
\end{equation}

Настоящая работа построена следующим образом. В
\S\,\ref{s2} формулируются
основные результаты работы.
Третий параграф  посвящен доказательству теоремы~\ref{t1},
четвертый~-- доказательству теоремы~\ref{t2}.
В~\S\,\ref{s5} приводятся некоторые
стандартные сведения
о~гиперэллиптических римановых поверхностях, используемые
в~работе, и выводится явное представление $\Psi$-функции
в~терминах абелевых дифференциалов.

В~заключение отметим, что близкие по содержанию результаты
получены в~работах \cite{19}, \cite{21},~\cite{32},
но при этом исходные данные в~них формулируются
непосредственно в~терминах чебышёвских коэффициентов $a_n$
и~$b_n$.

\goodbreak

\section{Формулировка основных результатов}
\label{s2}

\subsection{}
\label{s2.1}
Нам будет удобно перейти к~новым обозначениям. Будем
считать, что носитель меры~$\mu$ состоит из
непересекающихся отрезков $\Delta_j=[e_{2j-1},e_{2j}]$,
$j=\nobreak1,\dots,g+1$, расположенных на вещественной прямой
$\mathbb R$, $g\ge1$, $e_1<\dots<e_{2g+2}$, и конечного
числа точек $\zeta_k\in\mathbb R$, $k=1,\dots,m$. Меру~$\mu$
будем считать абсолютно непрерывной на
$E=\bigsqcup_{j=1}^{g+1}\Delta_j$ относительно меры Лебега
и такой, что
\begin{equation}
\label{eq18}
\frac{d\mu}{dx}=\frac1\pi\,
\frac{\rho(x)}{\sqrt{-h(x)}}>0 \quad \text{на } S,
\end{equation}
$E=\operatorname{ess\,supp}\mu$. В~\eqref{eq18} функция
$h(z)=\prod_{j=1}^{2g+2}(z-e_j)$, причем в~области
$D=\widehat{\mathbb C}\setminus S$ выбрана та ветвь
квадратного корня, которая положительна при положительных
значениях аргумента. Тем самым, $\sqrt{-h(x)}>0$ при
$x\in[e_{2g+1},e_{2g+2}]$; на остальных отрезках значение
корня определяется аналитическим продолжением. Весовая
функция $\rho$ предполагается \textit{голоморфной и отличной
от нуля на~$E$}. Вне~$E$ мера~$\mu$ имеет вид
\begin{equation}
\label{eq19}
d\mu(x)=\sum_{k=1}^m\lambda_k\delta(x-\zeta_k)\,dx,
\end{equation}
где все $\lambda_k>0$, точки
$\zeta_k\in\mathbb C\setminus\widehat E$,\enskip
$\widehat E$~-- выпуклая оболочка~$E$. Тем самым (см.~\eqref{eq2}),
\begin{equation}
\label{eq20}
\widehat\mu(z)=\frac1\pi\int_E\frac{\rho(x)}{z-x}\,\frac{dx}{\sqrt{-h(x)}}
+\sum_{k=1}^m\frac{\lambda_k}{z-\zeta_k}\,,
\qquad z\notin\operatorname{supp}\mu,
\end{equation}
$\operatorname{supp}\mu=E\cup\{\zeta_1,\dots,\zeta_m\}$.

Через $\omega_j(z)=\omega(z;\Delta_j,D)$, $j=1,\dots,g+1$,
условимся обозначать гармоническую меру (в~точке $z\in D$)
отрезка $\Delta_j$ относительно
области~$D=\widehat{\mathbb C}\setminus E$;
$g(z,\infty)=g_D(z,\infty)$~-- функция Грина области $D$
с~особенностью в~бесконечно удаленной точке $z=\infty$.

Пусть $\mathfrak R$ -- гиперэллиптическая риманова
поверхность рода $g$, заданная уравнением $w^2=h(z)$.
Будем считать, что $\mathfrak R$ реализована как
двулистное разветвленное в~точках $e_j$, $j=1,\dots,2g+2$,
накрытие римановой сферы $\widehat{\mathbb C}$ таким
образом, что переход с~одного листа на другой
осуществляется по верхнему~$\Delta^+_j$ и
нижнему~$\Delta^-_j$ берегам отрезков~$\Delta_j$. Тем самым, над
каждой точкой $\widehat{\mathbb C}$ за исключением точек
ветвления $e_j$ лежат ровно две точки римановой
поверхности, а~каждому отрезку $\Delta_j$ соответствует на
$\mathfrak R$ замкнутая аналитическая (в~комплексной
структуре $\mathfrak R$) кривая ${\boldsymbol\Gamma}_j$,
$j=1,\dots,g+1$, -- цикл на~$\mathfrak R$; положим
${\boldsymbol\Gamma}=\bigsqcup_{j=1}^{g+1}{\boldsymbol\Gamma}_j$.
Выбранная в~$D$ ветвь квадратного корня удовлетворяет
условию $\sqrt{h(z)}/z^{g+1}\to1$ при $z\to\infty$.
Функция $w$, $w^2=h(z)$, однозначна на $\mathfrak R$.
Первым (открытым) листом $D^{(1)}$ поверхности
$\mathfrak R$ будем считать тот, на котором
$w=\sqrt{h(z)}$. На втором листе $D^{(2)}$ имеем
$w=-\sqrt{h(z)}$. Для точек римановой
поверхности~$\mathfrak R$ будем использовать обозначение
$\mathbf z=(z,w)$, где $w=\pm\sqrt{h(z)}$; при этом для
точек первого листа $z^{(1)}=(z,\sqrt{h(z)}\,)$, а~для точек
второго $z^{(2)}=(z,-\sqrt{h(z)}\,)$. Вместо
$\mathbf z=(z,\pm\sqrt{h(z)}\,)$ иногда будем писать коротко
$\mathbf z=(z,\pm)$. Область $D^{(1)}$ будем, как правило,
отождествлять с~``физической'' областью $D$. Для
$\mathbf z=z^{(1)}$ будем иногда писать просто $w(z)$
вместо $w(\mathbf z)$; тем самым, приобретает смысл и
запись $w^\pm(x)=\sqrt{h(x\pm i0)}$, $x\in E$.
Каноническая проекция
$\operatorname{pr}\colon\mathfrak R\to\widehat{\mathbb C}$
определяется соотношением $\operatorname{pr}\mathbf z=z$,
в~частности,
$\operatorname{pr}D^{(1)}=\operatorname{pr}D^{(2)}=D$,
$\operatorname{pr}{\boldsymbol\Gamma}=E$. Замкнутые циклы
на $\mathfrak R$, соответствующие замкнутым лакунам
$[e_{2j},e_{2j+1}]$, $j=1,\dots,g$, будем обозначать через
$\mathbf L_j$. Тем самым,
$\operatorname{pr}\mathbf L_j=[e_{2j},e_{2j+1}]$.

Если на $\mathfrak R\setminus{{\boldsymbol\Gamma}}$ задана
функция $F(\mathbf z)$, то под $F^{(1)}(\mathbf x)$
понимаются ее ``некасательные'' предельные значения при
$z^{(1)}\to{\mathbf x}\in{\boldsymbol\Gamma}$,
$z^{(1)}\in D^{(1)}$, если они существуют; аналогичный
смысл придается и $F^{(2)}(\mathbf x)$:
$$
F^{(1)}(\mathbf x):=\lim_{z^{(1)}\to{\mathbf x}}F(z^{(1)}), \qquad
F^{(2)}(\mathbf x):=\lim_{z^{(2)}\to{\mathbf x}}F(z^{(2)}).
$$

Другие стандартные сведения о~гиперэллиптических римановых
поверхностях, которые нам здесь понадобятся, приведены
в~\cite{18} (приложения~A и~B), см.\
также~\cite{35},~\cite{36}.

\subsection{}
\label{s2.2}
Рассмотрим следующую систему из $g$ дифференциальных
уравнений относительно (неупорядоченного) набора $g$ точек
$\mathbf z_1(t),\dots,\mathbf z_g(t)$ на $\mathfrak R$,
$t\in\mathbb R$ (или, эквивалентно, относительно дивизора
$d(t)=\mathbf z_1(t)+\dots+\mathbf z_g(t)$ на
$S^g\mathfrak R$):
\begin{equation}
\label{eq21}
\frac{dz_j}{dt}=-\frac{w(\mathbf z_j)}{\prod_{k\ne j}(z_j-z_k)}
\int_{e_{2g+2}}^{+\infty}\frac{\prod_{k\ne j}(x-z_k)}{\sqrt{h(x)}}\,dx,
\qquad j=1,2,\dots,g.
\end{equation}
Как показано в~\cite{25} (теорема~2) (см.\
также \cite{18} (приложение~B, п.\,6)), если все
$z_j(0)=\operatorname{pr}\mathbf z_j(0)\in[e_{2j},e_{2j+1}]$,
$j=1,\dots,g$, то система~\eqref{eq21} интегрируется
в~явном виде преобразованием Абеля, ее решение задается
квазипериодической функцией непрерывного переменного~$t$
с~группой периодов $\omega_1(\infty),\dots,\omega_g(\infty)$
и со значениями в~$g$-мерном вещественном торе. Тем самым,
все $\mathbf z_j(t)\in\mathbf L_j$, $j=1,\dots,g$, при
$t>0$. Следовательно, все
$z_j(t)=\operatorname{pr}\mathbf z_j(t)\in[e_{2j},e_{2j+1}]$,
$j=1,\dots,g$, при $t>0$. Отметим, что система~\eqref{eq21}
приводит к~системе уравнений
Видома--Рахманова~\cite{37},~\cite{38}
\begin{gather}
\smash[b]{\sum_{j=1}^{g}\varepsilon_j\omega_k(z_j)
=(g-2t)\omega_k(\infty)-\frac2\pi\int_E\log|\rho(\zeta)|\,
\frac{\partial\omega_k(\zeta)}{\partial n^+_\zeta}\,d\zeta
-2\sum_{j=1}^m\omega_k(\zeta_j)\pmod2,}
\label{eq22}
\\[-1mm]
k=1,\dots,g,
\notag
\end{gather}
описывающих движение дивизора
$d(t)=\mathbf z_1(t)+\dots+\mathbf z_g(t)$ на торе
$\mathbf L_1\times\dotsb\times\mathbf L_g$, где
$\varepsilon_j=\pm1$ в~зависимости от точки
$\mathbf z_j=(z_j,\pm)$. Учитывая геометрический смысл
гармонической меры (см.~\cite{39},
а~также \cite{18}(приложение~B, п.\,7)),
уравнения~\eqref{eq22} естественно трактовать следующим
образом: дивизор
$d(t)=\mathbf z_1(t)+\dots+\mathbf z_g(t)$,
удовлетворяющий системе~\eqref{eq21}, движется с~постоянной
``угло\-вой'' скоростью (подробнее см.\
\cite{18}(приложение~B, п.\,6)).

Основными результатами настоящей работы являются
теоремы~\ref{t1}~и~\ref{t2}.

\begin{theorem}
\label{t1}
Пусть выполнены условия~\eqref{eq18}{\rm--}\eqref{eq19}.~Тогда
равномерно \text{внутри}
$\widehat{\mathbb C}\setminus\operatorname{supp}\mu$ имеем
\begin{equation}
\label{eq23}
G_{n+1}(z)=q_n(z)r_n(z)=\frac{\prod_{j=1}^g(z-z_j(n))}
{\sqrt{h(z)}}+o(\delta^n), \qquad n\to\infty,
\end{equation}
где $\delta\in(0,1)$ зависит от меры~$\mu$, величины
$z_1(n),\dots,z_g(n)$ соответствуют решению
$\mathbf z_1(t),\dots,\mathbf z_g(t)$ системы~\eqref{eq21},
взятому при $t=n\in\mathbb N$; начальные условия
$\mathbf z_1(0),\dots,\mathbf z_g(0)$ определяются мерой
$\mu$, при этом $z_j(0)\in[e_{2j},e_{2j+1}]$ для всех
$j=1,\dots,g$.
\end{theorem}

Доказательство теоремы~\ref{t1} основано на стандартной
технике~\cite{25}, \cite{20},~\cite{18}
исследования асимптотических свойств многочленов,
ортогональных на нескольких отрезках, и состоит в~сведении
задачи об асимптотике к~исследованию свойств
$\Psi$-функции и свойств решения (обобщенного)
\textit{сингулярного интегрального уравнения Наттолла}
(см.~\eqref{eq41}). Эти же свойства лежат в~основе
доказательства следующего результата.

\begin{theorem}
\label{t2}
Пусть выполнены условия~\eqref{eq18}{\rm--}\eqref{eq19}. Тогда
справедливы следующие асимптотические формулы следов для
коэффициентов $a_n$ и $b_n$ чебышёвской непрерывной
дроби~\eqref{eq3}:
\begin{gather}
b_n=\frac12\sum_{j=1}^{2g+2}e_j-\sum_{j=1}^gz_j(n-1)+o(\delta^n),
\qquad n\to\infty,
\label{eq24}
\\
\smash[b]{a_n=\operatorname{cap}{E}\cdot\exp\biggr\{\frac12\sum_{j=1}^g
g(\mathbf z_j(n),\infty) -\frac12\sum_{j=1}^g
g(\mathbf z_j(n-1),\infty)\biggl\}(1+o(\delta^n)),
\ n\to\infty,}
\label{eq25}
\end{gather}
где $\delta\in(0,1)$ зависит от меры~$\mu$.
\end{theorem}

Здесь $g(\mathbf z,\infty)$ -- однозначное продолжение
функции Грина $g_D(z,\infty)$ области~$D$ на всю риманову
поверхность~$\mathfrak R$, $\operatorname{cap}E$~--
логарифмическая емкость компакта~$E$.

Отметим, что при отсутствии у~меры~$\mu$ точечных масс
асимптотическая формула~\eqref{eq25} вытекает
непосредственно из результатов Видома~\cite{37}
(\S\,6, теорема~6.2 и \S\,9, теорема~9.1), но с~$o(1)$ вместо
$o(\delta^n)$ как в~\eqref{eq25}.

\subsection{}
\label{s2.3}
Рассмотрим на римановой поверхности~$\mathfrak R$
следующую краевую задачу Римана.

\begin{probl} 
\label{probR}
При фиксированном $n\in\mathbb N$, $n\ge g$, найти функцию
$\Psi=\Psi(\mathbf z;n)$ такую, что
\begin{enumerate}
\item 
\label{i1}
$\Psi$ (кусочно) мероморфна на
$\mathfrak R\setminus{\boldsymbol\Gamma}=D^{(1)}\sqcup{D^{(2)}}$;
\item 
\label{i2}
дивизор
$(\Psi)=(n-g)\infty^{(2)}+\mathbf z_1+\dots+\mathbf z_g+\zeta_1^{(1)}+\dots+\zeta_m^{(1)}-\zeta_1^{(2)}-\dots-\zeta_m^{(2)}-n\infty^{(1)}$;
\item 
\label{i3}
на ${\boldsymbol\Gamma}$ выполнено краевое
условие:
$\rho(x)\Psi^{(1)}(\mathbf x)=\Psi^{(2)}(\mathbf x)$,
$\mathbf x\in{\boldsymbol\Gamma}$.
\end{enumerate}
\end{probl}

В~п.\,\ref{i2} точки $\mathbf z_j$ -- ``свободные'' нули
$\Psi$-функции -- зависят от $n$, под
$\Psi^{(1)}(\mathbf x)$ в~п.\,\ref{i3} понимаются предельные
значения функции $\Psi(z^{(1)};n)$ при
$z^{(1)}\to\mathbf x\in{\boldsymbol\Gamma}$, аналогичный
смысл придается и $\Psi^{(2)}(\mathbf x)$. Так как
вес~$\rho$ марковский (см.~\eqref{eq18}) и все
$\lambda_k>0$, то $\mathbf z_j\in\mathbf L_j$, т.е.\
в~каждой лакуне $L_j=[e_{2j},e_{2j+1}]$ лежит ровно по одной
точке~$z_j$. При этом допускается, чтобы какая-нибудь из
точек $\mathbf z_j\in{\boldsymbol\Gamma}\cap\mathbf L_j$,
т.е.\ совпадала бы с~краем лакуны: $\mathbf z_j=e_{2j}$ или
$\mathbf z_j=e_{2j+1}$. В~таком случае эта точка считается
как нулем~$\Psi^{(1)}$, так и нулем~$\Psi^{(2)}$. При
этом, так как функция~$\rho$ голоморфна
на~${\boldsymbol\Gamma}$, сохраняются все нужные нам свойства
$\Psi$-функции, в~том числе
представление~\eqref{eq26}--\eqref{eq30} (см.\ ниже).

Функция $\Psi$, решающая задачу~\ref{probR}, всегда
существует. Так как род поверхности~$\mathfrak R$
положителен, то нули и полюсы~$\Psi$ на~$\mathfrak R$
оказываются связанными определенными соотношениями,
аналогичными соотношениям Абеля для мероморфной функции
на~$\mathfrak R$, а~дивизор $d=\mathbf z_1+\dots+\mathbf z_g$
является решением проблемы обращения Якоби. Анализ данных
этой проблемы показывает, что при
условиях~\eqref{eq18}--\eqref{eq19} на вес~$\rho$ и
постоянные~$\lambda_k$ ее решение всегда таково, что
$z_j=\operatorname{pr}{\mathbf z_j}\in[e_{2j},e_{2j+1}]$,
т.е.\ в~каждой замкнутой лакуне
$L_j=\operatorname{pr}{\mathbf L}$ между отрезками
$\Delta_1,\dots,\Delta_{g+1}$ лежит ровно по одной точке
$\operatorname{pr}{\mathbf z_j}$. Значит, дивизор
$d=\mathbf z_1+\dots+\mathbf z_g$ неспециальный,
а~следовательно, такая проблема обращения Якоби имеет
единственное решение. Отсюда вытекает, что $\Psi$-функция,
решающая задачу~\ref{probR}, единственна с~точностью до
нормировки и имеет в~бесконечно удаленной точке
$\mathbf z=\infty^{(1)}$ полюс в~точности $n$-го порядка.
При этом оказывается (подробнее см.~\S\,\ref{s3},
предложение~\ref{p1}), что равномерно внутри области
$\widehat{\mathbb C}\setminus({\widehat E}\cup\{\zeta_1,\dots,\zeta_m\})$,
где $\widehat E{\vrule width0pt height11pt}$~--   
выпуклая оболочка~$E$, справедлива
асимптотическая формула
$$
Q_n(z)=\Psi(z^{(1)};n)(1+o(1)), \qquad n\to\infty;
$$
здесь полиномы $Q_n(z)=\mathrm{const}\cdot q_n(z)$,
$\mathrm{const}\ne0$, нормированы так: старший коэффициент
$Q_n$ равен коэффициенту при степени~$z^n$ в~главной части
$\Psi$-функции, $z\to\infty$.

Нетрудно видеть, что для любого $z\notin E$ выполняется
соотношение
$$
\Psi(z^{(1)};n)\*\Psi(z^{(2)};n)
\equiv\mathrm{const}\prod_{j=1}^g(z-z_j),
\qquad
\mathrm{const}\ne0.
$$
В~дальнейшем мы будем
придерживаться следующей нормировки\footnote{Таким
условием $\Psi$-функция определяется
однозначно с~точностью до знака~``$\pm$'', выбор знака
уточняется в~дальнейшем.}
$\Psi$-функции:
$$
\Psi(z^{(1)};n)\Psi(z^{(2)};n)\equiv\prod_{j=1}^g(z-z_j).
$$
Нетрудно найти и явный вид этой функции. При
$\mathbf z\in\mathfrak R\setminus{\boldsymbol\Gamma}$ и
в~предположении, что $z_j\in(e_{2j},e_{2j+1})$, имеем:
\begin{equation}
\label{eq26}
\Psi(\mathbf z;n)=\Phi(\mathbf z)^{n-g}
e^{A(\mathbf z;\rho)}\mathscr F_n(\mathbf z)\Pi(\mathbf z).
\end{equation}
Здесь $\Phi(\mathbf z)=e^{G(\mathbf z,\infty)}$ --
(многозначная) отображающая функция,
$G(z,\infty)=g(z,\infty)+ig^*(z,\infty)$ -- комплексная
функция Грина области~$D$,
\begin{gather}
\label{eq27}
A(\mathbf z;\rho)=w(\mathbf z)\biggl\{\frac1{2\pi i}\int_E
\frac{\log\rho(x)}{z-x}\,\frac{dx}{w^{+}(x)}
+\frac12\,{c_{g+1}}
+\frac1{2\pi i}\sum_{k=1}^g
v_k\int_{\Delta_k}\frac1{z-x}\,\frac{dx}{w^{+}(x)}\biggr\},
\\[-4mm]
v_k=2\int_E\log\rho(x)\,d\Omega_k^+(x);
\nonumber
\end{gather}

\goodbreak

\noindent
при $g=0$
$$
\exp\{A(z^{(1)};\rho)\}=D(z;\rho)
=\exp\biggl(\frac{\sqrt{z^2-1}}{2\pi}\int_\Delta
\frac{\log\rho(x)}{x-z}\,\frac{dx}{\sqrt{1-\zeta^2}}\biggr)
$$
-- классическая функция Сегё. Функция
\begin{equation}
\label{eq28}
\mathscr F_n(\mathbf z)
=\exp\biggl(\,\sum_{j=1}^g\Omega(\mathbf z_j,\infty^{(1)};\mathbf z)
+2\pi i\sum_{k=1}^g\theta_k\Omega_k(\mathbf z)\biggr),
\end{equation}
величины
$\theta_k=\theta_k(n)=\ell_k(n)
+\boldsymbol\{(n-g)\omega_k(\infty)\boldsymbol\}$,
целые числа $\ell_k(n)\in\mathbb Z$ равномерно ограничены
при $n\to\infty$, а~дивизор
$d=\mathbf z_1+\dots+\mathbf z_g$, где
$\mathbf z_j=\mathbf z_j(n)$, является (единственным)
решением проблемы обращения Якоби
\begin{gather}
\label{eq29}
\smash[b]{\sum_{j=1}^g\Omega_k(\mathbf z_j)
\equiv\frac{i}\pi\int_E\log\rho(x)\,d\Omega^{+}_k(x)
-\sum_{j=1}^g\boldsymbol\{(n-g+\tfrac12)
\omega_j(\infty)\boldsymbol\}B_{kj}
+2\sum_{s=1}^m\Omega_k(\zeta_s^{(1)})},
\\[-1mm]
k=1,\dots,g.
\notag
\end{gather}
(Здесь и в~дальнейшем через
$\boldsymbol\{\cdot\boldsymbol\}$ обозначается дробная
часть соответствующего числа.) В~\eqref{eq27} величина
$c_{g+1}=c_{g+1}(n)$,
$e^{c_{g+1}}=\prod_{j=1}^g(e_{2g+2}-z_j)$. Отметим, что
равномерно ограниченные целые числа $\ell_k(n)$,
$k=1,\dots,g$, возникают в~\eqref{eq28} (см.\ формулу для
величин $\theta_k$) в~связи с~неоднозначностью абелевых
интегралов $\Omega_k(\mathbf z)$ для
$\mathbf z\in\partial\widetilde{\mathfrak R}$ при
интегрировании по путям, лежащим
в~$\widetilde{\mathfrak R}$,\enskip
$\widetilde{\mathfrak R}$ --
рассеченная риманова поверхность $\mathfrak R$ (подробнее
см.\ \cite{18} (приложения~A и~B)). Полином
$X_{g,n}(z):=\prod_{j=1}^g(z-z_j)$ является фактически
неизвестным ``полиномиальным параметром''
задачи~\ref{probR}. Наконец,
\begin{equation}
\label{eq30}
\Pi(\mathbf z)=\exp\biggl(\,\sum_{s=1}^m\Omega
\bigl(\zeta_s^{(1)},\zeta_s^{(2)};\mathbf z\bigr)\biggr).
\end{equation}
Так как функция $\rho$ голоморфна и отлична от нуля на
$E$, то правая часть представления~\eqref{eq26} имеет смысл
как голоморфная функция и при
$\mathbf z\in{D^{(1)}}\sqcup{\boldsymbol\Gamma}$ для
дивизора $d=\mathbf z_1+\dots+\mathbf z_g$,
удовлетворяющего условиям~\eqref{eq29}. Тем самым, под
функцией $\Psi(\mathbf z;n)$,
$\mathbf z\in{D^{(1)}}\sqcup{\boldsymbol\Gamma}$,
естественно понимать правую часть
представления~\eqref{eq26}. Аналогичное справедливо и для
функции $\Psi(\mathbf z;n)$ при
$\mathbf z\in{D^{(2)}}\sqcup{\boldsymbol\Gamma}$. На
${\boldsymbol\Gamma}$ эти два голоморфных продолжения не
совпадают: для них выполняется краевое условие п.\,\ref{i3}, где,
вообще говоря, функция $\rho\not\equiv1$.

Подробный вывод явных формул~\eqref{eq26}--\eqref{eq30} для
$\Psi$-функции при отсутствии точечных масс дан в
работе~\cite{18} (приложение~В). Для рассматриваемого
здесь общего случая вывод соответствующего представления
для $\Psi$-функции вполне аналогичен; ниже в~\S\,\ref{s5}
приводится краткая схема доказательства
представления~\eqref{eq26}--\eqref{eq30}.

\section{Доказательство теоремы~\ref{t1}}
\label{s3}

\subsection{}
\label{s3.1}
Для доказательства теоремы~\ref{t1} нам понадобится
рассмотреть и решить несколько более общую,
чем \ref{probR}, следующую ``модифицированную'' краевую
задачу \ref{probR'}. Пусть
$\zeta_j,\eta_j\in\mathbb R\setminus{\widehat E}$, $j=1,\dots,m$,
-- произвольные различные точки.

\goodbreak   

{\renewcommand{\theenumi}{\arabic{enumi}$'$)}
\renewcommand{\labelenumi}{\theenumi}
\begin{probll}
\label{probR'}
При фиксированном $n\in\mathbb N$, $n\ge g$, найти функцию
$\psi=\psi(\mathbf z;n)$ такую, что
\begin{enumerate}
\item 
\label{i1'}
$\psi$ (кусочно) мероморфна на
$\mathfrak R\setminus{\boldsymbol\Gamma}=D^{(1)}\sqcup{D^{(2)}}$;
\item 
\label{i2'}
дивизор
$(\psi)=(n-g)\infty^{(2)}+\mathbf z_1+\dots
+\mathbf z_g+\zeta^{(1)}_1+\dots+\zeta^{(1)}_m-y^{(2)}_1
-\dots\allowbreak-y^{(2)}_m-n\infty^{(1)}$;
\item 
\label{i3'}
на ${\boldsymbol\Gamma}$ выполнено краевое
условие:
$\rho(x)\psi^{(1)}(\mathbf x)=\psi^{(2)}(\mathbf x)$.
\end{enumerate}
\end{probll}}

Существование и единственность решения задачи
\ref{probR'} доказываются стандартным способом
(см.\ \cite{25},~\cite{17}, а~также \S\,\ref{s5} ниже):
из~предположения, что решение существует, вытекает его
явное представление; непосредственно проверяется, что полученная в
итоге функция является решением краевой задачи. Отметим,
что из условий этой задачи вытекает, что
\begin{equation}
\label{eq31}
\psi(z^{(1)};n)\psi(z^{(2)};n)
\equiv\mathrm{const}\cdot\frac{\prod_{j=1}^m(z-\zeta_j)}{\prod_{j=1}^m(z-\eta_j)}
\prod_{j=1}^g(z-z_j)
\end{equation}
при $z\in D$, где $\mathrm{const}\ne0$. Как обычно, будем
придерживаться той нормировки $\psi$-функции, которая дает
$\mathrm{const}=1$. Решение задачи~\ref{probR'}
понадобится нам при $n\to\infty$ и в~ситуации, когда
$\eta_j=\zeta_{j,n}\to \zeta_j$, $j=1,\dots,m$. Постановка
модифицированной задачи~\ref{probR'} связана с~наличием
у меры~$\mu$ точечных масс.

Соотношение~\eqref{eq15} определяет диагональную
аппроксимацию Паде $[n/n]_{\widehat\mu}=p_n/q_n$
функции~$\widehat\mu$. Непосредственно из теоремы
Рахманова~\cite{38} (\S\,2, теорема~1) (см.\
также~\cite{20}) о~сходимости диагональных
аппроксимаций Паде для рациональных возмущений марковских
функций вытекает, что в~условиях теоремы~\ref{t1} к~каждой
точке $\zeta_j$, $j=1,\dots,m$, при $n\to\infty$ стремится со
скоростью геометрической прогрессии ровно по одному полюсу
$n$-й диагональной аппроксимации Паде функции
$\widehat\mu$~-- нулю $\zeta_{j,n}$ ортогонального
полинома~$q_n$. \textit{Всюду в~дальнейшем предполагается,
что функция $\psi=\psi(\mathbf z;n)$~-- решение
задачи~\ref{probR'} при $n\ge g$ и
$\eta_j=\zeta_{j,n}\to \zeta_j$, $n\to\infty$}. Функция~$\psi$
определена однозначно и имеет следующий вид:
\begin{equation}
\label{eq32}
\psi(\mathbf z;n)=\Phi(\mathbf z)^{n-g}
e^{A(\mathbf z;\rho)}\mathscr F_n(\mathbf z)\Pi_n(\mathbf z).
\end{equation}
Здесь
\begin{gather}
\label{eq33}
A(\mathbf z;\rho)=w(\mathbf z)\biggl(\frac1{2\pi i}\int_E
\frac{\log\rho(x)}{z-x}\,\frac{dx}{w^{+}(x)}+\frac12\,{c_{g+1}}
+\frac1{2\pi i}\sum_{k=1}^g
v_k\int_{\Delta_k}\frac1{z-x}\, \frac{dx}{w^{+}(x)}\biggr),
\\[-5mm]
v_k=2\int_E\log\rho(x)\,d\Omega_k^+(x),
\nonumber
\\[1mm]
\label{eq34}
\mathscr F_n(\mathbf z)
=\exp\biggl(\,\sum_{j=1}^g\Omega(\mathbf z_j,\infty^{(1)};\mathbf z)
+2\pi i\sum_{k=1}^g\theta_k\Omega_k(\mathbf z)\biggr),
\end{gather}
где величины
$\theta_k=\theta_k(n)=\ell_k(n)+\boldsymbol\{(n-g)\omega_k(\infty)+\delta_n\boldsymbol\}$,
$\delta_n=o(\delta^n)$, $\delta\in(0,1)$, целые числа
$\ell_k(n)\in\mathbb Z$ равномерно ограничены при
$n\to\infty$, а~дивизор $d=\mathbf z_1+\dots+\mathbf z_g$
является (единственным)
решением проблемы обращения Якоби
\begin{align}
\sum_{j=1}^g\Omega_k(\mathbf z_j)
&\equiv\frac{i}\pi\int_E\log\rho(x)\,d\Omega^{+}_k(x)
-\sum_{j=1}^g\boldsymbol\{(n-g+\tfrac12)
\omega_j(\infty)\boldsymbol\}B_{kj}\notag
\\
&\qquad
+\sum_{s=1}^m\bigl(\Omega_k(\zeta_s^{(1)})+\Omega_k(\zeta_{s,n}^{(1)})\bigr),
\qquad k=1,\dots,g.
\label{eq35}
\end{align}
Наконец,
\begin{equation}
\label{eq36}
\Pi_n(\mathbf z)=
\exp\biggl(\,\sum_{s=1}^m\Omega\bigl(\zeta_s^{(1)},\zeta_{s,n}^{(2)};
\mathbf z\bigr)\biggr).
\end{equation}

Функция $[n/n]_{\widehat\mu}=p_n/q_n$ не зависит от
нормировки ортогональных многочленов. Имеем
$[n/n]_{\widehat\mu}=P_n/Q_n$, где полиномы $P_n,Q_n$,
$\deg{P_n}\le{n-1}$, $\deg{Q_n}\le n$, $Q_n\not\equiv0$,
могут быть найдены из соотношения
\begin{equation}
\label{eq37}
R_n(z):=(Q_n\widehat\mu-P_n)(z)=O\biggl(\frac1{z^{n+1}}\biggr),
\qquad z\to\infty.
\end{equation}
При таком подходе функцию $R_n$ принято называть
\textit{функцией остатка}, $Q_n$~-- знаменателем, а~$P_n$~--
числителем рациональной дроби $[n/n]_{\widehat\mu}$. Всюду
в~дальнейшем мы предполагаем, что $\psi$-функция задана
формулами~\eqref{eq32}--\eqref{eq36} (тем самым, выполняется
соотношение~\eqref{eq31} с~$\mathrm{const}=1$), а
полиномы~$Q_n=\mathrm{const}\cdot q_n$ \textit{нормированы
условием $Q_n(z)=\varkappa_nz^n+\dotsb$}, где
$\varkappa_n>0$~-- коэффициент при~$z^n$ в~разложении
функции $\psi(z^{(1)};n)$ в~ряд Лорана в~окрестности
бесконечно удаленной точки.

\subsection{}
\label{s3.2}
В~этом пункте мы выведем стандартным
образом~\cite{26},~\cite{20},~\cite{18} сингулярное
интегральное уравнение Наттолла для некоторой функции,
мероморфной на $\mathfrak R\setminus{{\boldsymbol\Gamma}}$
и связанной с~функциями $\psi$, $R_n$ и $Q_n$ (см.\ ниже
формулу~\eqref{eq39}), и покажем, что асимптотика свободных
нулей $R_n$ и ``блуждающих'' нулей $Q_n$ (иначе говоря,
``ложных'' полюсов диагональных аппроксимаций \text{Паде}
$P_n/Q_n$) полностью определяется асимптотическим при
$n\to\infty$ поведением точек
$\mathbf z_1(n),\dots,\mathbf z_g(n)$, удовлетворяющих
проблеме обращения Якоби~\eqref{eq35}, а~асимптотика самих
функций $R_n$ и $Q_n$ -- асимптотическим поведением
функции $\psi(\mathbf z;n)$.

Непосредственно из определения~\eqref{eq37} функции
остатка~$R_n$ и представления~\eqref{eq20} получаем, что
\begin{equation}
\label{eq38}
R^{+}_n(x)-R^{-}_n(x)=\frac{2Q_n(x)\rho(x)}{w^{+}(x)},
\qquad x\in E\setminus\{e_1,\dots,e_{2g+2}\}.
\end{equation}
Определим кусочно мероморфную на
$\mathfrak R\setminus{{\boldsymbol\Gamma}}$ функцию
$F(\mathbf z;n)$ следующим образом:
\begin{equation}
\label{eq39}
F(\mathbf z;n):=\psi(\mathbf z;n)\cdot\begin{cases}
R_n(z)w(z),&\mathbf z\in D^{(1)},
\\
2Q_n(z),&\mathbf z\in D^{(2)},
\end{cases}
\end{equation}
где функция $\psi(\mathbf z;n)$,
$\mathbf z\in{\mathfrak R\setminus{{\boldsymbol\Gamma}}}$,
определена в~\eqref{eq32}--\eqref{eq36}. Функция $F$ имеет
полюсы только в~точках $\infty^{(1)},\infty^{(2)}$, каждый
порядка $g$, а~из~\eqref{eq38} вытекает, что для скачка $F$
на ${\boldsymbol\Gamma}$ имеем
\begin{equation}
\label{eq40}
F^{(1)}(\mathbf x;n)-F^{(2)}(\mathbf x;n)
=V^{(2)}(\mathbf x)\frac1{\rho(x)}\,,
\qquad \mathbf x\in{\boldsymbol\Gamma},
\end{equation}
где функция
$V(\mathbf z)=-\psi(\mathbf z;n)R_n(z)w(z)
=\psi(\mathbf z;n)R_n(z)w(\mathbf z)$,
$\mathbf z\in D^{(2)}$, мероморфна на втором листе.
Действительно, положим $\psi_1(z):=\psi(z^{(1)};n)$,
$\psi_2(z):=\psi(z^{(2)};n)$, $z\,{\in}\,D$.  
Из краевого
условия~\ref{i3'} задачи~\ref{probR'} вытекает, что
$\rho(x)\psi^\pm_1(x)\,{=}\allowbreak\,\psi^\mp_2(x)$ 
при $x\in E$. Умножим
обе части~\eqref{eq38} на $w^+(x)\psi^+_1(x)$ и, пользуясь
последним краевым условием, преобразуем~\eqref{eq38} к~виду
$$
R^+_n(x)w^+(x)\psi^+_1(x)+\frac1{\rho(x)}R^-_n(x)w^-(x)\psi^-_2(x)
=2Q_n(x)\psi^-_2(x), \qquad x\in E.
$$
Аналогично получаем
$$
\frac1{\rho(x)}R^+_n(x)w^+(x)\psi^+_2(x)+R^-_n(x)w^-(x)\psi^-_1(x)
=2Q_n(x)\psi^+_2(x), \qquad x\in E.
$$
В~совокупности эти соотношения и дают краевые
условия~\eqref{eq40} для функции $F(\mathbf z;n)$,
заданной равенством~\eqref{eq39}.

Интегральная формула типа Коши~\cite{18} (формула~(A.11))
для $F$ принимает следующий вид:
\begin{equation}
\label{eq41}
F(\mathbf z;n)=-\frac1{2\pi i}\int_{{\boldsymbol\Gamma}^+}
V^{(2)}({\boldsymbol\zeta})\frac1{\rho(\zeta)}\,
d\Omega({\boldsymbol\zeta};\mathbf z)+p_n(z), \qquad
\mathbf z\notin{\boldsymbol\Gamma},
\quad \operatorname{deg}{p_n}\le{g},
\end{equation}
где дифференциал
$$
d\Omega({\boldsymbol\zeta};\mathbf z)
=\frac12\,\frac{w({\boldsymbol\zeta})+w(\mathbf z)}{\zeta-z}\,
\frac{d\zeta}{w({\boldsymbol\zeta})}\,.
$$
{\bad
Соотношение~\eqref{eq41} представляет собой
\textit{сингулярное интегральное уравнение Нат\-толла} для
функции~$F$ (ср.\ с~\cite{31}, где подобное уравнение
получено другим методом для классического случая $g=0$ и
$E=[-1,1]$).%
}
Рассмотрим эту формулу для
$\mathbf z\in{D^{(1)}}$. Так как функция $\rho(\zeta)$
голоморфна на каждой связной компоненте $\Delta_j$
компакта~$E$, то $\rho(\zeta)$ голоморфна на каждой кривой
${\boldsymbol\Gamma}_j$. Поэтому контур
${\boldsymbol\Gamma}$ можно, не меняя значения
интеграла~\eqref{eq41}, покомпонентно продеформировать
в~близкий контур $\Gamma^{(2)}$, также состоящий из
$(g+1)$-й компоненты и целиком расположенный на втором
листе~$D^{(2)}$ в~некоторой
окрестности~${\boldsymbol\Gamma}$, в~которой голоморфна и
отлична от нуля функция $\rho(\zeta)$ (напомним, что
функция $V({\boldsymbol\zeta})$ определена и мероморфна на
всем втором листе~$D^{(2)}$). Полученная интегральная
формула задает голоморфное продолжение функции
$F(\mathbf z;n)$, $\mathbf z\in D^{(1)}$, через контур
${\boldsymbol\Gamma}$ на второй лист римановой поверхности
вплоть до контура $\Gamma^{(2)}$ (отметим, что для
$\mathbf z\in D^{(2)}$ это продолжение не совпадает
с~функцией $F(\mathbf z;n)$, $\mathbf z\in D^{(2)}$,
определенной~\eqref{eq39}). Таким образом, получаем для
$z\in D$:
$$
F(z^{(1)};n)=-\frac1{2\pi i}
\int_{\Gamma^{(2)}}V(\zeta^{(2)})\frac1{\rho(\zeta)}\,
d\Omega(\zeta^{(2)};z^{(1)})
+p_n(z).
$$
Непосредственно из определения функций $F(\mathbf z;n)$ и
$V(\mathbf z)=V(z^{(2)})$ и с~учетом тождества
(см.~\eqref{eq31})
$$
\psi(z^{(1)};n)\psi(z^{(2)};n)
\equiv\frac{\prod_{j=1}^m(z-\zeta_j)}{\prod_{j=1}^m(z-\zeta_{j,n})}\,X_{g,n}(z)
$$
имеем для ${\boldsymbol\zeta}\in{D^{(2)}}$:
\begin{align*}
V({\boldsymbol\zeta})
&=V(\zeta^{(2)})=R_n(\zeta)w(\zeta^{(2)})\psi(\zeta^{(2)})
=-R_n(\zeta)w(\zeta^{(1)})\psi(\zeta^{(1)})
\frac{\psi(\zeta^{(2)})}{\psi(\zeta^{(1)})}
\\
&=-\psi(\zeta^{(1)})R_n(\zeta)w(\zeta^{(1)})
\frac{\psi(\zeta^{(1)})\psi(\zeta^{(2)})}
{\psi(\zeta^{(1)})^2}
=-F(\zeta^{(1)})\frac{T_n(\zeta)X_{g,n}(\zeta)}{\psi(\zeta^{(1)})^2}\,,
\end{align*}
где
$$
T_n(\zeta)\equiv\frac{\prod_{j=1}^m(\zeta-\zeta_j)}
{\prod_{j=1}^m(\zeta-\zeta_{j,n})}\to1,
\qquad n\to\infty,
$$
равномерно вне $D\setminus\{\zeta_1,\dots,\zeta_m\}$. Нетрудно
видеть, что операция инволюции $\mathbf z^*=(z,\mp w)$ при
$\mathbf z=(z,\pm w)$ обладает следующим свойством:
$d\Omega({\boldsymbol\zeta}^*,\mathbf z)
=d\Omega({\boldsymbol\zeta},\mathbf z^*)$.
Следовательно,
$d\Omega(\zeta^{(2)},z^{(1)})=d\Omega(\zeta^{(1)},z^{(2)})$,
и с~учетом сохранения положительной (относительно области
$D^{(1)}$) ориентации кривых интегрирования интегральное
представление~\eqref{eq41} для
$F_1(z):=F(z^{(1)};n)=\psi(z^{(1)};n)\*R_n(z)w(z^{(1)})$ и
$z\in D$ преобразуется к~следующему виду:
\begin{align}
F_1(z)
&=-\frac1{2\pi i}\int_{\Gamma^{(1)}}F(\zeta^{(1)})
\frac{T_n(\zeta)X_{g,n}(\zeta)}{\psi(\zeta^{(1)})^2}\,
\frac1{\rho(\zeta)}\,d\Omega(\zeta^{(1)};z^{(2)})+p_n(z)
\notag
\\
&=-\frac1{2\pi i}\int_{\Gamma^{(1)}}F_1(\zeta)
\frac{T_n(\zeta)X_{g,n}(\zeta)}{\psi(\zeta^{(1)})^2}\,
\frac1{\rho(\zeta)}\,d\Omega(\zeta^{(1)};z^{(2)})+p_n(z)
\notag
\\
&=-\frac1{2\pi i}\int_{{\boldsymbol\Gamma}^{+}}F^{(1)}(\zeta^{(1)})
\frac{T_n(\zeta)X_{g,n}(\zeta)}{[\psi^{(1)}(\zeta^{(1)})]^2}\,
\frac1{\rho(\zeta)}\,d\Omega^{+}(\zeta^{(1)};z^{(2)})+p_n(z).
\label{eq42}
\end{align}
Аналогично, для
$F_2(z):=F(z^{(2)};n)=2\psi(z^{(2)};n)Q_n(z)$ при $z\in D$
имеем:
\begin{equation}
\label{eq43}
F_2(z)
=-\frac1{2\pi i}\int_{{\boldsymbol\Gamma}^{+}}F^{(1)}(\zeta^{(1)})
\frac{T_n(\zeta)X_{g,n}(\zeta)}{[\psi^{(1)}(\zeta^{(1)})]^2}\,
\frac1{\rho(\zeta)}\,d\Omega^{+}(\zeta^{(1)};z^{(1)})+p_n(z)
\end{equation}
(ср.~\eqref{eq42}--\eqref{eq43} с~формулой~(4.2)
из~\cite{31}, где рассмотрен классический случай $g=0$,
$E=[-1,1]$). Подчеркнем, что эти представления справедливы
в~общем случае -- при любом расположении точек
$\mathbf z_j$ на~$\mathfrak R$. Действительно, при
$\mathbf z_j\ne{e_{2j},e_{2j+1}}$ функции, стоящие под
знаком интеграла в~\eqref{eq42}--\eqref{eq43}, голоморфны
на ${\boldsymbol\Gamma}_j$ и ${\boldsymbol\Gamma}_{j+1}$,
так как функция $\psi$ не обращается в~нуль на этих
компонентах ${\boldsymbol\Gamma}$. Если же $\mathbf z_j$
совпадает с~краем лакуны, то нуль в~знаменателе,
происходящий от функции $\psi^2$, компенсируется нулем
функции $F(\mathbf z;n)$ на ${\boldsymbol\Gamma}$ в~точке
$\mathbf z_j$ и нулем полинома $X_{g,n}(z)$.

Формулы~\eqref{eq42}--\eqref{eq43} лежат в~основе наших
последующих рассуждений. Поясним, какого рода информацию
мы собираемся из них извлечь. Предположим, например, что
по некоторой подпоследовательности
$\Lambda\subset\mathbb N$ имеем:
$z_j(n)\Subset(e_{2j},e_{2j+1})$ при $n\to\infty$,
$j=1,\dots,g$. Тогда из~\eqref{eq42} мы получим, что для
$n\in\Lambda$ асимптотики функции $F_1(z):=F(z^{(1)};n)$ и
полинома $p_n$ в~области $D$ совпадают. Функция $F_1(z)$
имеет в~$D$ нули в~тех точках $z_j$, для которых
$\mathbf z_j=z_j^{(1)}$. Следовательно, полином $p_n$
обращается в~нуль в~близких точках. Аналогичный результат
вытекает из~\eqref{eq43} для полинома $p_n$ и функции
$F_2(z):=F(z^{(2)};n)$. Из определения функций $F_1$ и
$F_2$ вытекает, что их нули, порожденные нулями
$\mathbf z_1,\dots,\mathbf z_g$ функции $\psi$, различны
между собой. Следовательно, для достаточно больших
$n\in\Lambda$ полином $p_n$, степень которого $\le g$,
имеет ровно $g$ нулей в~$D$. Тем самым,
$\operatorname{deg}{p_n}=g$ и при должной нормировке
асимптотика $p_n$ совпадает с~асимптотикой полинома
$\prod_{j=1}^g(z-z_j(n))$. Отсюда уже вытекает асимптотика
оставшихся нулей функций $F_1$ и $F_2$ в~$D$. Но их нулями
в~$D$, отличными от точек $z_j$, могут быть только нули
$R_n$ и $Q_n$ соответственно. Тем самым,
$$
Q_n(z)R_n(z)
=\varkappa_n^2\biggl(\frac{\prod_{j=1}^g(z-z_j(n))}
{\sqrt{h(z)}}+o(1)\biggr), \qquad n\to\infty, \quad n\in\Lambda.
$$
Следовательно, равномерно внутри
$D\setminus\{\zeta_1,\dots,\zeta_m\}$
$$
q_n(z)r_n(z)
=\frac{\prod_{j=1}^g(z-z_j(n))}
{\sqrt{h(z)}}+o(1), \qquad n\to\infty, \quad n\in\Lambda.
$$

\subsection{}
\label{s3.3}
Изложим теперь приведенные выше соображения более
формально и в~полной общности, не ограничиваясь каким-либо
частным случаем в~поведении точек $z_j(n)$. Подчеркнем,
что цель наших исследований -- доказать, что
асимптотическое поведение нулей
$\mathbf z_1(n),\dots,\mathbf z_g(n)$ функции $\psi$
\textit{полностью определяет} асимптотическое поведение
других возможных нулей функции $F$, которыми, как следует
из~\eqref{eq39}, могут быть лишь нули функции остатка~$R_n$
или ортогонального полинома~$Q_n$.

Для произвольного $\varepsilon>0$
положим
$D_\varepsilon=\{z\in\widehat{\mathbb C}:
g(z,\infty)>\varepsilon\}$,
$E_\varepsilon=\{z\in\nobreak\widehat{\mathbb C}:
g(z,\infty)<\varepsilon\}$
-- окрестность компакта~$E$,
$\Gamma_\varepsilon=\partial{D_\varepsilon}=\{z\in\widehat{\mathbb C}:g(z,\infty)=\varepsilon\}$,
$D_0=D$. Выберем $\varepsilon_0>0$ так, чтобы все
$\zeta_j\in D_{\varepsilon_0}$, а~функция $\rho\ne0$ и
голоморфна в~$E_{\varepsilon_0}$. В~дальнейшем мы
рассматриваем только $\varepsilon\in(0,\varepsilon_0)$.
Зафиксируем произвольное $\varepsilon\in(0,\varepsilon_0)$
и некоторую последовательность $\{\varepsilon_n\}$ такую,
что $\varepsilon_n\in(\varepsilon/3,2\varepsilon/3)$ и при
$n\to\infty$
\begin{equation}
\label{eq44}
\max_{z\in\Gamma_{\varepsilon_n}}|p_n(z)|
=O\Bigl(\,\min_{z\in\Gamma_{\varepsilon_n}}|p_n(z)|\Bigr),
\qquad
\max_{z\in\Gamma_{\varepsilon_n}}|X_{g,n}(z)|
=O\Bigl(\,\min_{z\in\Gamma_{\varepsilon_n}}|X_{g,n}(z)|\Bigr).
\end{equation}
Тем самым, $X_{g,n}(z)=\prod_{j=1}^g(z-z_j)\ne0$ и
$p_n(z)\ne0$ при $z\in\Gamma_{\varepsilon_n}$. Так как
функция $F_1(z):=F(z^{(1)};n)$ голоморфна
в~$D\setminus\{\infty\}$ и дифференциал
$d\Omega(\zeta^{(1)};z^{(2)})$ не имеет особенности при
$z,\zeta\in D$ (см.\ \cite{18} (формула~(A.12))), то
в~\eqref{eq42} для $z\in D$ мы можем заменить контур
интегрирования ${\boldsymbol\Gamma}\subset\mathfrak R$ на
контур $\Gamma_{\varepsilon_n}$:
\begin{equation}
\label{eq45}
F_1(z)
=-\frac1{2\pi i}\int_{\Gamma_{\varepsilon_n}}F_1(\zeta)
\frac{T_n(\zeta)X_{g,n}(\zeta)}{\psi(\zeta^{(1)})^2}\,
\frac1{\rho(\zeta)}\,d\Omega(\zeta^{(1)};z^{(2)})+p_n(z),
\end{equation}
причем формула~\eqref{eq45} будет справедлива и для
$z\in\Gamma_{\varepsilon_n}=\partial{D_{\varepsilon_n}}$.
В~таком случае, учитывая найденное ранее представление для
функции~$\psi$ и выбор величин~$\varepsilon_n$
(см.~\eqref{eq44}), получаем, что равномерно по
$z\in\Gamma_{\varepsilon_n}$ выполняется следующее
соотношение:
\begin{equation}
\label{eq46}
|F_1(z)-p_n(z)|
=o(1)\cdot\max_{\zeta\in\Gamma_{\varepsilon_n}}|F_1(\zeta)|
=o(1)\cdot M_n,
\end{equation}
где
$$
M_n=\max_{\zeta\in\Gamma_{\varepsilon_n}}|F_1(\zeta)|,
\qquad o(1)=O(\delta^n), \quad
\delta=\delta_n=e^{-2\varepsilon_n}<e^{-2\varepsilon/3}<1.
$$
Нетрудно видеть, что из~\eqref{eq46} вытекает соотношение
$$
m_n:=\max_{z\in\Gamma_{\varepsilon_n}}|p_n(z)|=M_n(1+o(1)).
$$
Тем самым, из~\eqref{eq46} с~учетом~\eqref{eq44} получаем
$$
|F_1(z)-p_n(z)|
=o(1)\cdot m_n=o(1)\min_{z\in\Gamma_{\varepsilon_n}}|p_n(z)|<|p_n(z)|,
\qquad z\in\Gamma_{\varepsilon_n}.
$$
Отметим, что непосредственно отсюда по теореме Руше
вытекает, что для функций $F_1$ и $p_n$ разность числа их
нулей и полюсов в~области $D_{\varepsilon_n}$ одинакова.
Пусть теперь $t=3\varepsilon/4>\varepsilon_n$,
$$
m_n(t)=\max_{z\in\Gamma_t}|p_n(z)|, \qquad
M_n(t)=\max_{z\in\Gamma_t}|F_1(z)|.
$$
Так как $p_n$ -- полином степени $\le g$, то
$m_n\le m_n(t)$. Аналитическая функция
$(F_1(z)-p_n(z))/\Phi^g(z)$ многозначна в~$D$,
но имеет однозначный модуль (здесь и далее
$\Phi(z)=e^{G(z,\infty)}$ -- отображающая функция).
Поэтому к~ней применим принцип максимума модуля.
Следовательно,
$$
\max_{z\in\Gamma_t}\biggl|\frac{F_1(z)-p_n(z)}{\Phi^g(z)}\biggr|
\le\max_{z\in\Gamma_{\varepsilon_n}}
\biggl|\frac{F_1(z)-p_n(z)}{\Phi^g(z)}\biggr|
=O(1)\max_{z\in\Gamma_{\varepsilon_n}}|F_1(z)-p_n(z)|.
$$
Тем самым,
\begin{equation}
\label{eq47}
\max_{z\in\Gamma_{t}}|F_1(z)-p_n(z)|
=o(1)\cdot m_n=o(1)\cdot m_n(t).
\end{equation}
Отсюда уже легко вытекает, что
$$
M_n(t)=m_n(t)\cdot(1+o(1)).
$$

Заменим теперь в~\eqref{eq43} для $F_2(z):=F(z^{(2)};n)$
контур интегрирования
${\boldsymbol\Gamma}\subset\mathfrak R$ на контур
$\Gamma_{\varepsilon_n}\subset D$ и рассмотрим эту формулу
для $z\in\Gamma_t$ (напомним, что
$t=3\varepsilon/4>2\varepsilon/3>\varepsilon_n$):
$$
F_2(z)
=-\frac1{2\pi i}\int_{\Gamma_{\varepsilon_n}}F_1(\zeta)
\frac{T_n(\zeta)X_{g,n}(\zeta)}{\psi(\zeta^{(1)})^2}\,
\frac1{\rho(\zeta)}\,d\Omega(\zeta^{(1)};z^{(1)})+p_n(z),
\qquad z\in\Gamma_t.
$$
Тогда, учитывая выбор параметра $t$ и полученные выше
соотношения между $M_n$, $m_n$, $M_n(t)$ и $m_n(t)$, получаем
\begin{equation}
\label{eq48}
\max_{z\in\Gamma_{t}}|F_2(z)-p_n(z)|
=o(1)\cdot M_n=o(1)\cdot m_n=o(1)\cdot m_n(t).
\end{equation}
Введем теперь временно \textit{новую нормировку\/}: положим
$m_n(t)=1$ и сохраним прежние обозначения для остальных
величин. Тогда из~\eqref{eq46} и~\eqref{eq48} для
многозначных аналитических в~$D$ функций
$F_1(z)/\Phi^g(z)$, $F_2(z)/\Phi^g(z)$ и
$p_n(z)/\Phi^g(z)$, имеющих в~$D$ однозначные модули,
получаем
\begin{align}
&\biggl|\frac{F_1(z)}{\Phi^g(z)}-\frac{p_n(z)}{\Phi^g(z)}\biggr|
=o(1) \quad \text{равномерно по $z\in\overline{D_t}$},
\label{eq49}
\\[2mm]
&\biggl|\frac{F_2(z)}{\Phi^g(z)}-\frac{p_n(z)}{\Phi^g(z)}\biggr|
=o(1) \quad \text{равномерно по $z\in\overline{D_t}$}
\label{eq50}
\end{align}
(предполагается, что в~левых частях~\eqref{eq49}
и~\eqref{eq50} выбирается одна и та же ветвь многозначной
функции~$\Phi$). Рассмотрим последовательность функций
$\{p_n/\Phi^g\}$. По принципу максимума модуля для
аналитических функций с~однозначным модулем имеем:
$|p_n(z)/\Phi^g(z)|\le e^{-tg}$ при $z\in\overline{D_t}$
равномерно по $n$ (здесь мы учли новую нормировку
$m_n(t)=1$). С~другой стороны, применяя принцип максимума
модуля к~полиномам $p_n$ в~области $E_\tau$ при
произвольном $\tau>t$, получаем
$\max_{z\in\Gamma_\tau}|p_n(z)|>m_n(t)=1$. Тем самым,
функция, тождественно равная нулю, не является предельной
точкой последовательности $\{p_n/\Phi^g\}$ (в~топологии
равномерной сходимости на компактных подмножествах области~$D_t$).
Следовательно, если мы применим
к~последовательностям функций $\{F_1/\Phi^g\}$,
$\{F_2/\Phi^g\}$ и $\{p_n/\Phi^g\}$ теорему Гурвица, то
в~силу соотношений~\eqref{eq49} и~\eqref{eq50} получим, что
асимптотическое поведение нулей этих функций,
расположенных в~области~$D_t$, одинаково. Так как
$|\Phi(z)|\ne0$ в~$D$ и любая однозначная в~окрестности
точки $z=\infty$ ветвь $\Phi$ имеет полюс первого порядка
в~точке $z=\infty$, то, значит, асимптотическое поведение
нулей и полюсов функций $F_1$, $F_2$ и $p_n$,
расположенных в~области~$D_t$, также одинаково. Поскольку
$\operatorname{deg}{p_n}\le g$, то отсюда
сразу следует, что функции $F_1$ и $F_2$ \textit{имеют
в~$D_t$ не более чем по $g$ нулей}.

Сделаем теперь дальнейшие выводы из того, что функции
$F_1$ и $F_2$ имеют асимптотически одинаковое поведение
нулей в~$D_t$. Из явного вида~\eqref{eq39} функции
$F(\mathbf z;n)$ вытекает, что возможные нули функции
$F_1$ -- это или точки $z_j$ при $\mathbf z_j=z_j^{(1)}$,
или нули функции остатка $R_n$, а~возможные нули функции
$F_2$ -- это или точки $z_j$ при $\mathbf z_j=z_j^{(2)}$,
или нули полинома $Q_n$. В~силу доказанного и~те,
и~другие точки порождают асимптотически близкие к~ним нули
полинома $p_n$. Если общее число асимптотически различных
нулей функций $F_1$ и $F_2$ равно~$g$, то, значит, имеется
ровно $g$ близких к~ним нулей полинома $p_n$. Так как
$\operatorname{deg}{p_n}\le g$, то получаем полное
асимптотическое описание полинома $p_n$,
$\operatorname{deg}{p_n}=g$, в~терминах нулей $F_1$ и $F_2$.

Предположим теперь, что для некоторой последовательности
$\Lambda\subset\mathbb N$ имеем: все
$z_j(n)\Subset(e_{2j},e_{2j+1})$, $j=1,\dots,g$,
$n\in\Lambda$. Так как $t=3\varepsilon/4$, а
$\varepsilon>0$ произвольно, то для достаточно малого
$\varepsilon$ все $z_j\in{D_t}$ и попарно различны между
собой. Те точки $z_j$, для которых
$\mathbf z_j=z^{(1)}_j$, являются нулями $F_1$ и в~силу
доказанного ``притягивают'' нули $p_n$. Остальные $z_j$,
$\mathbf z_j=z^{(2)}_j$, -- это нули $F_2$, которые также
``притягивают'' нули $p_n$. Так как все точки
$z_1,\dots,z_g$ различны между собой и число их равно $g$,
то мы получили полное описание всех нулей $p_n$. Точнее,
$\operatorname{deg}{p_n}=g$ для всех достаточно больших
$n\in\Lambda$ и $\alpha^{-1}_np_n(z)=X_{g,n}(z)+o(1)$,
$n\in\Lambda$, $n\to\infty$, где $\alpha_n$ -- старший
коэффициент $p_n$.

Наконец, рассмотрим случай, когда для некоторого
$k\in\{1,\dots,g\}$ точка $z_k(n)$ стремится к~краю
лакуны: $z_k(n)\to{e_{2k}}$ или $z_k(n)\to{e_{2k+1}}$ при
$n\to\infty$. Покажем, что и в~этом случае
$\operatorname{deg}{p_n}=g$ для достаточно больших~$n$,
а~нули~$p_n$ бесконечно близки к~точкам
$z_1(n),\dots,z_g(n)$ при $n\to\infty$. Выберем, как и выше,
в~представлении~\eqref{eq45} для функции $F_1(z)$
в~качестве~$\Gamma^{(1)}$ кривую
$\Gamma^{(1)}_\varepsilon=\{z:g(z,\infty)=\varepsilon\}$,
$\varepsilon\in(0,\varepsilon_0)$, и запишем это
представление для $z\in D$ в~следующем виде
(см.~\eqref{eq45}):
$$
F_1(z)
=-\frac1{2\pi i}\int_{\Gamma^{(1)}_\varepsilon}F_1(\zeta)
\frac{T_n(\zeta)X_{g,n}(\zeta)}{\psi(\zeta^{(1)})^2}\,
\frac1{\rho(\zeta)}\,d\Omega(\zeta^{(1)};z^{(2)})
+p_n(z).
$$
Будем считать, что этой формулой функция $F_1$ задана на
$D^{(1)}$, т.е.\ для $\mathbf z\in{D^{(1)}}$:
\begin{equation}
\label{eq51}
F_1(\mathbf z)=-\frac1{2\pi i}\int_{\Gamma^{(1)}_{\varepsilon}}
F_1(\zeta^{(1)})
\frac{T_n(\zeta)X_{g,n}(\zeta)}{\psi(\zeta^{(1)})^2}\,
\frac1{\rho(\zeta)}\,d\Omega(\zeta^{(1)};\mathbf z^{*})+p_n(z).
\end{equation}
В~силу свойств дифференциала
$d\Omega({\boldsymbol\zeta};\mathbf z)$ правая часть
этого равенства голоморфно продолжается по $\mathbf z$ на
$D^{(2)}$ вплоть до кривой $\Gamma^{(2)}_\varepsilon$,
лежащей на втором листе $D^{(2)}$ ``над'' кривой
$\Gamma_\varepsilon=\Gamma^{(1)}_\varepsilon$. Тем самым,
функция $F_1(\mathbf z)$, $\mathbf z\in{D^{(1)}}$,
\textit{голоморфно продолжается на второй лист~$D^{(2)}$
вплоть до кривой $\Gamma^{(2)}_\varepsilon$} (это
продолжение, вообще говоря, не совпадает с~функцией
$F(\mathbf z;n)$, $\mathbf z\in D^{(2)}$). С~учетом
определения функции $F_1(z):=F(z^{(1)};n)$ получаем, что
на второй лист голоморфно продолжается произведение
$\psi(z^{(1)};n)R_n(z)w(z)$, $z\in D$. Обратимся теперь
к~краевым условиям на функцию $\psi$:
$\rho(\mathbf x)\psi^{(1)}(\mathbf x)=\psi^{(2)}(\mathbf x)$,
$\mathbf x\in{\boldsymbol\Gamma}$. Функция~$\rho$
голоморфна и отлична от нуля в~некоторой окрестности~$E$.
Поэтому, заменив в~этих краевых условиях кривую
${\boldsymbol\Gamma}$ на кривую $\Gamma^{(2)}_\varepsilon$
при произвольном $\varepsilon\in(0,\varepsilon_0)$, мы,
как легко увидеть, получим краевую задачу, решением
которой является та же самая функция~$\psi$ (это видно и
непосредственно из явных формул~\eqref{eq32}--\eqref{eq36}).
Таким образом, функция $\psi(z^{(1)};n)$, $z\in D$,
голоморфно продолжается на второй лист вплоть до кривой
$\Gamma^{(2)}_\varepsilon$. Выберем теперь
$\varepsilon_n\in(\varepsilon/3,2\varepsilon/3)$ так,
чтобы при $n\to\infty$ выполнялись
соотношения~\eqref{eq44}. Рассуждения, вполне аналогичные
приведенным выше и опирающиеся на принцип максимума модуля
и теорему Гурвица для многозначных аналитических функций
в~области
$D(n)=D^{(1)}\sqcup{\boldsymbol\Gamma}\sqcup S^{(2)}_{\varepsilon_n}$,
показывают, что в~этой ситуации имеет место аналог
соотношения~\eqref{eq49}:
\begin{equation}
\label{eq52}
\biggl|\frac{F_1(\mathbf z)}{\Phi^g(\mathbf z)}
-\frac{p_n(\mathbf z)}{\Phi^g(\mathbf z)}\biggr|
=o(1) \quad \text{равномерно по }
\mathbf z\in\overline{D(n)},
\end{equation}
а все нули полинома $p_n$ ``притягиваются'' к~точкам
$z_1,\dots,z_n$ при $n\to\infty$. Поясним более подробно
ту часть рассуждений, которая касается аналитического
продолжения функций с~первого листа $D^{(1)}$ через кривую
${\boldsymbol\Gamma}$ на второй лист~$D^{(2)}$.
Предположим для определенности, что $z_k\to{e_{2k}}$ при
$n\to\infty$, тем самым, фактически речь идет
о~продолжении через кривую ${\boldsymbol\Gamma}_k$,
$\operatorname{pr}{{\boldsymbol\Gamma}_k}=\Delta_k=[e_{2k-1},e_{2k}]$.
Будем считать для простоты, что $\Delta_k=[-1,1]$. Тогда
с~помощью функции $Z=z+\sqrt{z^2-1}$, обратной функции
Жуковского, двулистная гиперэллиптическая риманова
поверхность $\mathfrak R$ конформно отображается на
риманову сферу $\widehat{\mathbb C}$, из которой удалены
$2g$ отрезков вида $[\alpha_j,\beta_j]$,
$[1/\beta_j,1/\alpha_j]$, $j=1,\dots,g$, где все отрезки
$[\alpha_j,\beta_j]$ лежат во внешности единичного круга
$|Z|>1$, а~соответствующие им парные отрезки
$[1/\beta_j,1/\alpha_j]$ -- внутри него. При этом верхние
берега парных отрезков отождествляются (``склеиваются'')
между собой; аналогичное правило действует и для нижних
берегов. При таком отождествлении получаем, очевидно,
сферу с~$g$ ручками. Первому (открытому) листу $D^{(1)}$
римановой поверхности $\mathfrak R$ соответствует
внешность единичного круга $|Z|>1$ с~разрезами по отрезкам
$[\alpha_j,\beta_j]$, $j=1,\dots,g$, второму листу
$D^{(2)}$ -- внутренность единичного круга $|Z|<1$ с
разрезами по отрезкам $[1/\beta_j,1/\alpha_j]$,
$j=1,\dots,g$; см.\ \cite{18} (рис.~1). Обратное
отображение с~помощью функции Жуковского
$z=\frac12(Z+1/Z)$ задает параметрическое представление
$\mathfrak R$ и делает очевидным процесс аналитического
продолжения через кривую
${\boldsymbol\Gamma}_k=\{z(Z):|Z|=1\}$ с~первого листа
римановой поверхности на второй. Действительно, для
весовой функции $\rho$ имеем
$\rho(z)=\rho\bigl(\frac12(Z+1/Z)\bigr)$, для полинома
$p_n$ имеем $p_n(z)=p_n\bigl(\frac12(Z+1/Z)\bigr)$. Отметим,
что из такого симметричного представления вытекает, что
каждой точке $z_k\to{e_{2k}}$ соответствуют два близких
нуля продолженной функции $p_n\bigl(\frac12(Z+1/Z)\bigr)$.
Это означает, что продолженная на второй лист функция
$F_1(\mathbf z)$ также имеет в~окрестности точки $e_{2k}$
два нуля: один нуль происходит от функции
$\psi(\mathbf z;n)$, а~второй -- от функции остатка
$R_n(z)$ (см.~\eqref{eq39}).

Итак, мы показали, что в~общей ситуации, во-первых,
$\operatorname{deg}{p_n}=g$ для достаточно больших $n$ и,
во-вторых, асимптотика $g$ нулей полинома $p_n$ полностью
определяется асимптотическим поведением точек
$z_1(n),\dots,z_g(n)$, точнее,
$$
\alpha^{-1}_np_n(z)=X_{g,n}(z)+o(1), \qquad n\to\infty,
$$
где $\alpha_n$ -- старший коэффициент $p_n$. С~учетом
последнего равенства нетрудно видеть, что полученные выше
соотношения~\eqref{eq49}--\eqref{eq50} сохранятся, если
вместо нормировки $m_n(t)=1$ мы будем считать, что старший
коэффициент полинома~$p_n$ равен~$2$:
$p_n(z)=2z^g+\dotsb$. Тем самым, при $n\to\infty$
\begin{equation}
\label{eq53}
\begin{alignedat}{2}
p_n(z)&=2X_{g,n}(z)+o(1)
&\qquad &\text{равномерно внутри } \mathbb C,
\\
p_n(z)&=2X_{g,n}(z)(1+o(1)) &\qquad &\text{равномерно
на } \Gamma_{\varepsilon_n}.
\end{alignedat}
\end{equation}
Следовательно, во всех полученных выше соотношениях $p_n$
можно заменить на $2X_{g,n}$. Кроме того, как нетрудно
видеть, в~этих соотношениях величина $o(1)=o(\delta^n)$,
где $\delta\in(0,1)$ и зависит от~$\mu$.

\subsection{}
\label{s3.4}
С~помощью соотношений~\eqref{eq49}--\eqref{eq50}
и~\eqref{eq53} уже легко завершается доказательство
теоремы~\ref{t1}. Действительно, в~силу определения
функций $F_1$ и $F_2$ из~\eqref{eq49}--\eqref{eq50}
и~\eqref{eq53} получаем, что при $n\to\infty$ равномерно
на $\Gamma_{\varepsilon_n}$
\begin{equation}
\label{eq54}
\begin{aligned}
\psi(z^{(1)};n)R_n(z)w(z)&=2X_{g,n}(z)(1+o(1)),
\\
\psi(z^{(2)};n)Q_n(z)&=X_{g,n}(z)(1+o(1)).
\end{aligned}
\end{equation}
Перемножим теперь левые и правые части этих соотношений и
воспользуемся тем, что
$\psi(z^{(1)};n)\psi(z^{(2)};n)=T_n(z)X_{g,n}(z)$, где
$T_n(z)\to1$ для $z\notin{\zeta_j}$. Тогда получим
$$
T_n(z)Q_n(z)R_n(z)
=\frac{2X_{g,n}(z)}{w(z)}(1+o(1))
=\frac{2X_{g,n}(z)}{w(z)}+o(1)
\quad \text{равномерно на } \Gamma_{\varepsilon_n}.
$$
Но $\varepsilon_n\in(\varepsilon/3,2\varepsilon/3)$, тем
самым, равномерно
в~$\overline{D_\varepsilon}=\{z:g(z,\infty)\ge\varepsilon\}$
имеем
\begin{equation}
\label{eq55}
T_n(z)Q_n(z)R_n(z)=\frac{2X_{g,n}(z)}{w(z)}+o(1),
\qquad n\to\infty.
\end{equation}
Так как
$$
Q_n(z)R_n(z)
=\frac1\pi\int_E\frac{Q_n^2(x)}{z-x}\,\frac{\rho(x)\,dx}{\sqrt{-h(x)}}
+\sum_{k=1}^m\frac{\lambda_kQ_n^2(\zeta_k)}{z-\zeta_k}
$$
и $T_n(z)=1+o(1)$, то из~\eqref{eq55} вытекает, что
\begin{equation}
\label{eq56}
\frac1\pi\int_E\frac{Q_n^2(x)\rho(x)\,dx}{\sqrt{-h(x)}}
+\sum_{k=1}^m\lambda_kQ_n^2(\zeta_k)=2+o(1), \qquad n\to\infty.
\end{equation}
Тем самым, окончательно получаем
$$
T_n(z)q_n(z)r_n(z)
=\frac{X_{g,n}(z)}{\sqrt{h(z)}}+o(\delta^n), \qquad n\to\infty, \quad \delta\in(0,1),
$$
равномерно в
$\overline{D_\varepsilon}=\{z:g(z,\infty)\ge\varepsilon\}$.
Поскольку $\varepsilon>0$ произвольно, а
$T_n(z)=1+o(1)$ равномерно вне
$\widehat{\mathbb C}\setminus\{\zeta_1,\dots,\zeta_m\}$, то
$$
q_n(z)r_n(z)
=\frac{X_{g,n}(z)}{\sqrt{h(z)}}+o(\delta^n), \qquad n\to\infty,
$$
равномерно внутри $D\setminus\{\zeta_1,\dots,\zeta_m\}$.
Теорема~\ref{t1} доказана.

\subsection{}
\label{s3.5}
Отметим, что при подходящей нормировке полиномов $Q_n(z)$~--
знаменателей диагональных аппроксимаций Паде функции
$\widehat\mu$~-- их сильная асимптотика внутри области
$\widehat{\mathbb C}\setminus\widehat E$ и на компакте~$E$
описывается в~терминах $\Psi$-функции, решающей
задачу~\ref{probR}, следующим образом.

\begin{prop}
\label{p1}
Пусть голоморфная на~$E$ функция $\rho$ удовлетворяет
ус\-ло\-вию~\eqref{eq18}, $\Psi=\Psi(\mathbf z;n)$ -- решение
задачи~\ref{probR} при $n\in\mathbb N$,
$\Psi_1(z):=\Psi(z^{(1)};n)$. Тогда при подходящей
нормировке полиномов $Q_n(z)$ имеем:
\begin{enumerate}
\item 
\label{i1p}
$Q_n(z)\,{=}\,\Psi_1(z)(1+o(\delta^n))$,  
$n\,{\to}\,\infty$,   
равномерно внутри
$\widehat{\mathbb C}\setminus(\widehat E\cup\{\zeta_1,\dots,\zeta_m\})$;
\item 
\label{i2p}
$Q_n(x)X_{g,n}(x)=(\Psi_1^+(x)+\Psi_1^-(x))X_{g,n}(x)+o(\delta^n)$,
$n\to\infty$, равномерно на~$E$.
\end{enumerate}
\end{prop}

Здесь под $\Psi_1^+(x)$ ($\Psi_1^-(x)$) понимаются верхние
(соответственно нижние) предельные значения $\Psi_1(z)$
на~$E$. Функция $\Psi$, решающая задачу~\ref{probR},
имеет в~точке $\mathbf z=\infty^{(1)}$ полюс в~точности
$n$-го порядка. Так как $\operatorname{deg}{Q_n}=n$ при
всех~$n$ и
$\Psi(z^{(1)};n)/z^n\to K_n\in\mathbb R\setminus\{0\}$ при
$z\to\infty$, то в~условиях предложения~\ref{p1}
естественно нормировать $Q_n$ условием
$Q_n(z)=K_nz^n+\dotsb$. С~помощью
представления~\eqref{eq26}--\eqref{eq30} величину~$K_n$
нетрудно найти и в~явном виде
(см.\ \cite{18} (приложение~B)), при этом, так как
$K_n\in\mathbb R\setminus\{0\}$, можно считать, что
$K_n>0$; таким выбором знака у~$K_n$ однозначно
определяется и сама $\Psi$-функция.

\begin{proof}[предложения~\ref{p1}]
Пункт~\ref{i1p} предложения вытекает непосредственно из
соотношения~\eqref{eq54} и тождества
$$
\Psi(z^{(1)};n)\Psi(z^{(2)};n)=X_{g,n}(z).
$$

Для доказательства п.\,\ref{i2p} поступим следующим образом.
Из~\eqref{eq52} и~\eqref{eq53} вытекает, что равномерно на
${\boldsymbol\Gamma}$ имеем
$$
F_1(\mathbf z)=2X_{g,n}(z)+o(1),\qquad n\to\infty,
\quad
o(1)=o(\delta^n).
$$
Пользуясь определением функций $F_1$ и~$F$,
перепишем последнее соотношение в~виде двух равенств
для $x\in E$ следующим образом:
$$
\Psi_1^+(x)R^+_n(x)w^+(x)=2X_{g,n}(x)+o(1),
\qquad
\Psi_1^-(x)R^-_n(x)w^-(x)=2X_{g,n}(x)+o(1).
$$
Умножим обе части первого равенства на $\Psi_1^-(x)$,
второго -- на $\Psi_1^+(x)$, сложим получившиеся
соотношения и воспользуемся равенствами~\eqref{eq38} и
$\rho(x)\*\Psi_1^+(x)\*\Psi_1^-(x)=X_{g,n}(x)$. Получим
$$
Q_n(x)X_{g,n}(x)=\bigl(\Psi_1^+(x)+\Psi_1^-(x)\bigr)X_{g,n}(x)+o(1),
$$
что и требовалось доказать.
\end{proof}

\section{Доказательство теоремы~\ref{t2}}
\label{s4}

\subsection{}
\label{s4.1}
Первая часть теоремы~\ref{t2}~-- формула~\eqref{eq24}~--
вытекает непосредственно из теоремы~\ref{t1}.
Действительно, с~одной стороны,
\begin{equation}
\label{eq57}
q_n(z)r_n(z)
=\frac1\pi\int_E\frac{q_n^2(x)}{z-x}\,\frac{\rho(x)\,dx}{\sqrt{-h(x)}}
+\sum_{k=1}^m\frac{\lambda_kq_n^2(\zeta_k)}{z-\zeta_k}\,,
\end{equation}
с~другой --
\begin{equation}
\label{eq58}
q_n(z)r_n(z)
=\frac{X_{g,n}(z)}{\sqrt{h(z)}}+o(\delta^n),
\qquad n\to\infty, \quad \delta\in(0,1),
\end{equation}
равномерно внутри $D'_\varepsilon$, где
$D'_\varepsilon=\{z:g(z,\infty)>\varepsilon\}\setminus\{\zeta_1,\dots,\zeta_m\}$,
$\varepsilon>0$ произвольно. Пользуясь~\eqref{eq57},
разложим функцию $q_n(z)r_n(z)$ в~ряд по степеням $1/z$
в~бесконечно удаленной точке:
$$
q_n(z)r_n(z)
=\frac{d_0}z+\frac{d_1}{z^2}+\dotsb,
$$
где
\begin{align*}
d_0&=\frac1\pi\int_E\frac{q_n^2(x)\rho(x)\,dx}{\sqrt{-h(x)}}
+\sum_{k=1}^m\lambda_kq_n^2(\zeta_k)=1,
\\
d_1&=\frac1\pi\int_E\frac{xq_n^2(x)\rho(x)\,dx}{\sqrt{-h(x)}}
+\sum_{k=1}^m\lambda_k\zeta_kq_n^2(\zeta_k).
\end{align*}
Непосредственно из рекуррентной формулы~\eqref{eq11}
вытекает, что $d_1=b_{n+1}$. В~силу~\eqref{eq58}
коэффициенты $d_0,d_1,\dots$ с~точностью до $o(\delta^n)$
совпадают с~соответствующими лорановскими коэффициентами
$c_0,c_1,\dots$ функции $X_{g,n}(z)/\sqrt{h(z)}$. Так как
$c_0=1$, то
$$
c_1=\lim_{z\to\infty}
z\biggl(\frac{zX_{g,n}(z)}{\sqrt{h(z)}}-1\biggr)
=\lim_{z\to\infty}
z\biggl(\frac{zX_{g,n}(z)-\sqrt{h(z)}}{z^g}\,\biggr).
$$
Прямые вычисления дают для этой величины следующее
представление:
$$
c_1=\frac12\sum_{j=1}^{2g+2}e_j-\sum_{j=1}^gz_j(n).
$$
Тем самым,
$$
b_{n+1}=\frac12\sum_{j=1}^{2g+2}e_j-\sum_{j=1}^gz_j(n)+o(\delta^n),
\qquad n\to\infty.
$$
Формула~\eqref{eq24} доказана.

\subsection{}
\label{s4.2}
Для доказательства~\eqref{eq25} поступим следующим образом.
Из рекуррентной формулы~\eqref{eq11} вытекает, что
\begin{equation}
\label{eq59}
a_n=\frac1\pi\int_E
xq_{n-1}(x)q_n(x)\frac{\rho(x)\,dx}{\sqrt{-h(x)}}
+\sum_{k=1}^m\lambda_k\zeta_kq_{n-1}(\zeta_k)q_n(\zeta_k)
=\frac{k_{n-1}}{k_n}\,,
\end{equation}
$k_n>0$ -- старший коэффициент ортонормированного полинома
$q_n$. Из~\eqref{eq56} получаем, что $k_n$ и
$\varkappa_n>0$ -- старший коэффициент полинома $Q_n$ --
связаны соотношением $\varkappa_n=k_n(\sqrt2+o(1))$, где
$o(1)=o(\delta^n)$. Следовательно,
$a_n=\varkappa_{n-1}/\varkappa_n(1+o(1))$. Но величина
$\varkappa_n=\varkappa_n(\psi)$ связана с~нормировкой
$\psi$-функции: $\psi(z^{(1)};n)/z^n\to\varkappa_n$ при
$z\to\infty$. Пользуясь явным
представлением~\eqref{eq32}--\eqref{eq36} для $\psi$-функции,
найдем теперь асимптотику отношения величин
$\varkappa_{n-1}$ и $\varkappa_n$ при $n\to\infty$.

Так как $\zeta_{j,n}\to \zeta_j$ со скоростью геометрической
прогрессии, то (см.\ п.\,\ref{s3.5})
$\varkappa_n(\psi)=K_n(\Psi)(1+o(\delta^n))$. Таким
образом, достаточно найти асимптотику отношения величин
$K_{n-1}$ и~$K_n$ при $n\to\infty$.

В~силу тождества
$$
\Psi(z^{(1)};n)\Psi(z^{(2)};n)\equiv\prod_{j=1}^g(z-z_j)
$$
имеем
$$
\Psi(z^{(2)};n)=\frac1{K_nz^{n-g}}+\dotsb, \qquad z\to\infty.
$$
Следовательно,
$$
\frac{\Psi(z^{(1)};n)}{\Psi(z^{(2)};n)}
=K^2_nz^{2n-g}+\dotsb, \qquad z\to\infty,
$$
и
$$
K^2_n=\lim_{z\to\infty}
\frac{\Psi(z^{(1)};n)}{\Psi(z^{(2)};n)}\,z^{g-2n},
$$
при этом можно считать, что $z>0$. С~помощью явных
формул~\eqref{eq26}--\eqref{eq29} получаем
\begin{equation}
\label{eq60}
\frac{\Psi(z^{(1)};n)}{\Psi(z^{(2)};n)}
=\Phi^{2n-2g}e^{A(z^{(1)};\rho)-A(z^{(2)};\rho)}
\frac{\mathscr F_n(z^{(1)})}{\mathscr F_n(z^{(2)})}\,,
\end{equation}
где первый сомножитель ведет себя как
$(\operatorname{cap}E)^{2g-2n}z^{2n-2g}+\nobreak\dotsb$, при $z\to\infty$
а~второй
не зависит от $n$. Тем самым,
\begin{equation}
\label{eq61}
\frac{K^2_{n-1}}{K^2_n}
=\lim_{z\to\infty}\biggl[\frac{z^{2}}{\Phi^{2}(z)}\cdot\frac{\mathscr F_{n-1}(z^{(1)})}{\mathscr F_{n-1}(z^{(2)})}\cdot\frac{\mathscr F_{n}(z^{(2)})}{\mathscr F_{n}(z^{(1)})}\biggr].
\end{equation}
Преобразуем теперь величину
${\mathscr F_n(z^{(1)})}/{\mathscr F_n(z^{(2)})}
=\exp(\varphi_n(z))$,
где (см.~\eqref{eq29})
\begin{align}
\varphi_n(z)
&:=\sum_{j=1}^g
\Bigl[\Omega\bigl(\mathbf z_j(n),\infty^{(1)};z^{(1)}\bigr)
-\Omega\bigl(\mathbf z_j(n),\infty^{(1)};z^{(2)}\bigr)\Bigr]
\notag
\\
&\qquad+2\pi i\sum_{j=1}^g
\theta_j(n)\bigl(\Omega_j(z^{(1)}) -\Omega_j(z^{(2)})\bigr).
\label{eq62}
\end{align}
По правилу перестановки пределов интегрирования и
параметров для $a$-нор\-ми\-ро\-ван\-ных абелевых дифференциалов
(см.\ \cite{18} (формула~(A.7))) имеем
\begin{align*}
&\Omega(\mathbf z_j,\infty^{(1)};z^{(1)})
-\Omega(\mathbf z_j,\infty^{(1)};z^{(2)})
=-\int_{z^{(1)}}^{z^{(2)}}d\Omega(\mathbf z_j,\infty^{(1)};{\boldsymbol\zeta})
\\
&\qquad
=-\int^{\infty^{(1)}}_{\mathbf z_j}d\Omega(z^{(1)},z^{(2)};{\boldsymbol\zeta})
\;(\operatorname{mod}{2\pi i})
=\int_{\infty^{(1)}}^{\mathbf z_j}d\Omega(z^{(1)},z^{(2)};{\boldsymbol\zeta})
\;(\operatorname{mod}{2\pi i}).
\end{align*}
Так как при $x\in\mathbb R\setminus E$
$$
dG(\zeta^{(1)},\zeta^{(2)},{\boldsymbol\zeta})
+2\pi i\sum_{j=1}^g\omega_j(x)\,d\Omega_j({\boldsymbol\zeta})
+d\Omega(\zeta^{(1)},\zeta^{(2)},{\boldsymbol\zeta})=0
$$
(см.\ \cite{18} (формула~(A.6))), то при $z>e_{2g+2}$
\begin{align}
\int_{\infty^{(1)}}^{\mathbf z_j}d\Omega(z^{(1)},z^{(2)};{\boldsymbol\zeta})
={}&{-}\int_{\infty^{(1)}}^{\mathbf z_j}
dG(z^{(1)},z^{(2)};{\boldsymbol\zeta})
-2\pi i\sum_{k=1}^{g}\omega_k(z)
\int_{\infty^{(1)}}^{\mathbf z_j}d\Omega_k({\boldsymbol\zeta})
\notag
\\
={}&{-}\bigl(g(\mathbf z_j,z)-g(\infty^{(1)},z)\bigr)+i\beta_j(n)
\notag  
\\
&{-}2\pi i\sum_{k=1}^g\omega_k(z)
\bigl(\Omega_k(\mathbf z_j)-\Omega_k(\infty^{(1)})\bigr),
\label{eq63}
\end{align}
где $\beta_j(n)\in\mathbb R$. Воспользуемся теперь тем,
что
$\Omega_k(z^{(1)})-\Omega_k(z^{(2)})=2\Omega_k(z^{(1)})$,
и~\cite{18} (формула~(B.12)). Тогда из~\eqref{eq62}
и~\eqref{eq63} получим
\begin{align}
\varphi_n(z)
={}&{-}\sum_{j=1}^{g}\bigl[g(\mathbf z_j,z)-g(\infty^{(1)},z)\bigr]
+2\pi i\sum_{k=1}^g\omega_k(z)
\biggl(\,\sum_{j=1}^g\theta_j(n)B_{kj}\biggr)
\notag
\\
&+2\pi i\sum_{j=1}^g\theta_j(n)\cdot2\Omega_j(z^{(1)})+i\beta(n)+C_1(z),
\label{eq64}
\end{align}
где $\beta(n)\in\mathbb R$,
$\theta_j=\theta_j(n)=\ell_j(n)+\boldsymbol\{(n-g)\omega_j(\infty)\boldsymbol\}$,
величины $\ell_j(n)\in\mathbb Z$ и равномерно ограничены
при $n\to\infty$, а~$C_1(z)$ не зависит от~$n$. Поменяем
в~\eqref{eq64} порядок суммирования во~втором слагаемом и
воспользуемся тем, что (см.\ \cite{18} (формула~(A.8)))
$$
\sum_{k=1}^g\omega_k(z)B_{kj}
=-2i\operatorname{Im}\Omega_j(z^{(1)}),
\qquad z\in D.
$$
Тогда это слагаемое примет вид
\begin{align*}
&2\pi i\sum_{k=1}^g\omega_k(z)\biggl(\,\sum_{j=1}^g\theta_jB_{kj}\biggr)
=2\pi i\sum_{j=1}^g\theta_j\biggl(\,\sum_{k=1}^g\omega_k(z)B_{kj}\biggr)
\\
&\qquad\qquad
=2\pi i\sum_{j=1}^g\theta_j
\bigl(-2i\operatorname{Im}\Omega_j(z^{(1)})\bigr)
=-2\pi i\sum_{j=1}^g\theta_j
\bigl(2i\operatorname{Im}\Omega_j(z^{(1)})\bigr).
\end{align*}
Следовательно, из~\eqref{eq64} имеем
$$
\varphi_n(z)=-\sum_{j=1}^gg(\mathbf z_j,z)+i\beta(n)+C_2(z),
$$
где $\beta(n)\in\mathbb R$, а~$C_2(z)$ не зависит от $n$.
Окончательно из последней формулы и~\eqref{eq61} получаем
$$
\frac{K^2_{n-1}}{K^2_n}
=(\operatorname{cap}E)^2\exp\biggl(\,\sum_{j=1}^gg
\bigl(\mathbf z_j(n),\infty\bigr)
-\sum_{j=1}^gg\bigl(\mathbf z_j(n-1),\infty\bigr)\biggr),
$$
тем самым,
\begin{align*}
a_n &=\frac{K_{n-1}}{K_n}(1+o(\delta^n))
\\
&=\operatorname{cap}{E}\cdot\exp\biggl(\frac12\sum_{j=1}^g
g\bigl(\mathbf z_j(n),\infty\bigr) -\frac12\sum_{j=1}^g
g\bigl(\mathbf z_j(n-1),\infty\bigr)\biggr)(1+o(\delta^n)).
\end{align*}
Асимптотическая формула~\eqref{eq25}, а~вместе с~ней и
теорема~\ref{t2} доказаны.

\section({Доказательство существования, единственности и
вывод явных формул для \$\000\134psi\$-функции})
{Доказательство существования, единственности и
вывод явных формул для $\psi$-функции}
\label{s5}

\subsection*{Единственность}
Для доказательства единственности решения краевой
задачи~\ref{probR'}
предположим, что существуют два
решения этой задачи: $\psi$ с~дивизором свободных нулей
$d=\mathbf z_1+\dots+\mathbf z_g$ и $\widetilde\psi$
с~дивизором свободных нулей
$\widetilde d=\widetilde{\mathbf z}_1+\dots+\widetilde{\mathbf z}_g$.
Так как $\mathbf z_j,\widetilde{\mathbf z}_j\in\mathbf L_j$,
то каждый из дивизоров $d$ и $\widetilde d$
неспециальный. Функция
$F(\mathbf z;n):=\psi(\mathbf z;n)/\widetilde\psi(\mathbf z;n)$~--
однозначная мероморфная функция на всей римановой
поверхности~$\mathfrak R$ с~дивизором
$(F)=d-\widetilde d=\mathbf z_1+\dots
+\mathbf z_g-\widetilde{\mathbf z}_1-\dots-\widetilde{\mathbf z}_g$.
Из хорошо известного свойства
(см.\ \cite{35}, \cite{18} (приложение~A, п.\,4))
эквивалентных дивизоров, заданных на гиперэллиптической
римановой поверхности, вытекает, что
$\mathbf z_j=\widetilde{\mathbf z}_j$, $j=1,2,\dots,g$, тем
самым, $(F)=0$ и $F\equiv\mathrm{const}$.

\goodbreak

\subsection*{Существование}
Доказательство существования $\psi$-функции при условии,
что $\eta_j=\zeta_{j,n}$ и $\zeta_{j,n}\to \zeta_j$ при $n\to\infty$,
проведем по следующей стандартной схеме. Предположив
сначала, что все $z_j\in(e_{2j},e_{2j+1})$, найдем явный
вид этой функции в~терминах абелевых интегралов на~$\mathfrak R$.
После чего можно непосредственно проверить,
что найденная функция является решением
задачи~\ref{probR'}, а~явная формула дает искомое
решение задачи и в~случае, когда хотя бы одна точка~$z_j$
совпадает с~краем лакуны $e_{2j}$ или~$e_{2j+1}$. При этом
окажется, что точки $\mathbf z_1,\dots,\mathbf z_g$
удовлетворяют проблеме обращения
Якоби \cite{18} (формула~(A.2)) с~правой частью,
имеющей специальный вид (см.~\eqref{eq35}).
В~\cite{18} (приложение~B) показано,
что такая специальная
задача всегда имеет единственное решение
и~притом такое, что
все $z_j\in[e_{2j},e_{2j+1}]$.

\subsection{}
\label{s5.1}
Итак, пусть $\psi$ -- решение задачи~\ref{probR'},
причем такое, что все $z_j\in(e_{2j},e_{2j+1})$. Положим
$v(z):=\psi(z^{(1)};n)\psi(z^{(2)};n)$,
$z\in D'=D\setminus\{\zeta_1,\dots,\zeta_m\}$. Легко увидеть,
что эта функция продолжается на всю риманову сферу~$\widehat{\mathbb C}$
как (однозначная) мероморфная функция
с~полюсом порядка~$g$ в~точке $z=\nobreak\infty$, простыми
полюсами в~точках~$\zeta_{1,n},\dots,\zeta_{m,n}$ и простыми
нулями в~точках~$\zeta_1,\dots,\zeta_m$ и точках $z_1,\dots,z_g$.
Следовательно, $v(z)$~-- рациональная функция, точнее,
$$
v(z)\equiv\mathrm{const}\cdot
\frac{\prod_{k=1}^m(z-\zeta_k)}{\prod_{k=1}^m(z-\zeta_{k,n})}
\prod_{j=1}^g(z-z_g),
$$
где $\mathrm{const}\ne0$. Для последующих рассуждений
удобно считать, что $\mathrm{const}=1$, тем самым,
\begin{equation}
\label{eq65}
\psi(z^{(1)};n)\psi(z^{(2)};n)\equiv
\frac{\prod_{k=1}^m(z-\zeta_k)}{\prod_{k=1}^m(z-\zeta_{k,n})}
\prod_{j=1}^g(z-z_g).
\end{equation}
Этим соотношением $\psi$-функция определена однозначно
с~точностью до знака~``$\pm$''; в~дальнейшем мы уточним
выбор знака. Отметим, что фактически
$X_{g,n}(z)=\prod_{j=1}^g(z-z_j)$ является
``полиномиальным параметром'' задачи~\ref{probR'}. Для
функции $\psi_1(z):=\psi(z^{(1)};n)$ на
$E=\bigsqcup_{j=1}^{g+1}\Delta_j$ выполняется следующее
краевое условие:\footnote{В~классическом
случае $g=0$ непосредственно из
этого краевого условия уже вытекает нужное представление
функции~$\psi$.}
\begin{equation}
\label{eq66}
\rho(x)\psi_1^{+}(x)\psi_1^{-}(x)
=\frac{\prod_{k=1}^m(z-\zeta_k)}{\prod_{k=1}^m(z-\zeta_{k,n})}
\prod_{j=1}^g(x-z_j), \qquad x\in{E}.
\end{equation}
Действительно, пусть $z\to x\in\Delta_j$,
$j\in\{1,\dots,g+1\}$, причем $\operatorname{Im}z>0$.
Тогда $\psi(z^{(1)};n)\to\psi^+_1(x)$ и в~силу краевого
условия
$\rho(x)\psi^{(1)}(\mathbf x)=\psi^{(2)}(\mathbf x)$,
$\mathbf x\in{\boldsymbol\Gamma}$, имеем
$\psi(z^{(2)};n)\to\rho(x)\psi^-_1(x)$. Теперь
\eqref{eq66}~вытекает непосредственно из~\eqref{eq65}.

\subsection{}
\label{s5.2}
Стандартным образом \cite{18} (приложение~B, пп.\,1--4)
устанавливается следующее представление для функции
$\psi(z^{(1)};n)$, справедливое в~предположении, что все
$z_j\in(e_{2j},e_{2j+1})$:
{\allowdisplaybreaks
\begin{align}
\psi(z^{(1)};n)
&=e^{(n-g)G(z,\infty)}
\exp\biggl(-\sum_{j=1}^g\Omega(\infty^{(1)},\mathbf z_j;z^{(1)})
+\sum_{k=1}^m\Omega\bigl(\zeta_k^{(1)},\zeta_{k,n}^{(2)};z^{(1)}\bigr)\biggr)\notag
\\*
&\qquad\times
\exp\biggl(\frac{w(z)}{2\pi i}
\int_E\frac{\log\rho(x)\,dx}{(z-x)w^{+}(x)}
+\frac{c_{g+1}}2+v(z)w(z)\biggr)
\notag
\\
&\qquad\times
\exp\biggl(2\pi i\sum_{j=1}^g
\bigl((n-g)\omega_j(\infty)+m_j\bigr)\Omega_j(z^{(1)})\biggr)
\notag
\\
&=\Phi^{n-g}(z)
\exp\biggl(\,\sum_{j=1}^g\Omega(\mathbf z_j,\infty^{(1)};z^{(1)})
+\sum_{k=1}^m\Omega\bigl(\zeta_k^{(1)},\zeta_{k,n}^{(2)};z^{(1)}\bigr)\biggr)
\notag
\\
&\qquad\times
\exp\biggl(\frac{w(z)}{2\pi i}\int_E\frac{\log\rho(x)\,dx}{(z-x)
w^{+}(x)}+\frac{c_{g+1}}2+v(z)w(z)\biggr)
\notag
\\
&\qquad\times
\exp\biggl(2\pi i\sum_{j=1}^g\theta_j(n)\Omega_j(z^{(1)})\biggr),
\label{eq67}
\end{align}}
где $\Phi(z)=e^{G(z,\infty)}$ -- (многозначная)
отображающая функция,
$$
\theta_k=\theta_k(n)
=\ell_k(n)+\boldsymbol\{(n-g)\omega_k(\infty)+\delta_n\boldsymbol\},
\qquad
\delta_n=o(\delta^n),
\quad
\delta\in(0,1),
$$
целые числа
$\ell_k(n)\in\mathbb Z$ равномерно ограничены при
$n\to\infty$,
$$
v(z):=
\frac1{2\pi i}\sum_{k=1}^gv_k\int_{\Delta_k}\frac1{z-x}\,
\frac{dx}{w^+(x)}\,,
\qquad
v_k=2\int_E\log\rho(x)\,d\Omega_k^+(x).
$$

Нетрудно видеть, что для функции
$\psi_2(z):=\psi(z^{(2)};n)$, $z\in D$, выполняется
краевое условие
$\psi_2^+(x)\psi_2^-(x)=\rho(x)T_n(x)X_{g,n}(x)$, $x\in E$
(ср.\ с~\cite{18} (формула~(B.2))).
Отсюда уже легко
вытекает следующее представление для $\psi$-функ\-ции при
$\mathbf z=\nobreak z^{(2)}$:
{\allowdisplaybreaks
\begin{align}
\psi(z^{(2)};n)
&=e^{(g-n)G(z,\infty)}
\exp\biggl(-\sum_{j=1}^g\Omega(\infty^{(1)},\mathbf z_j;z^{(2)})
+\sum_{k=1}^m\Omega\bigl(\zeta_k^{(1)},\zeta_{k,n}^{(2)};z^{(2)}\bigr)\biggr)
\notag
\\*
&\qquad\times
\exp\biggl(\frac{w(z)}{2\pi i}\int_E\frac{\log\rho(x)\,dx}{(x-z)w^{+}(x)}
-\frac{c_{g+1}}2-v(z)w(z)\biggr)
\notag
\\*
&\qquad\times
\exp\biggl(2\pi i\sum_{j=1}^g
\bigl((n-g)\omega_j(\infty)+m_j\bigr)\Omega_j(z^{(2)})\biggr)
\notag
\\
&=\frac1{\Phi^{n-g}(z)}
\exp\biggl(\,\sum_{j=1}^g\Omega(\mathbf z_j,\infty^{(1)};z^{(2)})
+\sum_{k=1}^m\Omega\bigl(\zeta_k^{(1)},\zeta_{k,n}^{(2)};z^{(2)}\bigr)\biggr)
\notag
\\*
&\qquad\times
\exp\biggl(\frac{w(z)}{2\pi i}\int_E\frac{\log\rho(x)\,dx}{(x-z)w^{+}(x)}
-\frac{c_{g+1}}2-v(z)w(z)\biggr)
\notag
\\*
&\qquad\times
\exp\biggl(2\pi i\sum_{j=1}^g\theta_j(n)\Omega_j(z^{(2)})\biggr).
\label{eq68}
\end{align} } 

Подведем итог. Мы показали, что в~предположении
$z_j\in(e_{2j},e_{2j+1})$, $j=1,\dots,g$, существование
решения задачи~\ref{probR'} эквивалентно тому, что
точки $\mathbf z_1,\dots,\mathbf z_g$ являются решением
проблемы обращения Якоби~\eqref{eq35} со специальной правой
частью. Одновременно мы нашли явное
представление~\eqref{eq67}--\eqref{eq68} функции
$\psi(\mathbf z;n)$ для
$\mathbf z\in\mathfrak R\setminus{\boldsymbol\Gamma}$ и
при нормировке
$$
\psi(z^{(1)};n)\psi(z^{(2)};n)
=\frac{\prod_{k=1}^m(z-\zeta_k)}{\prod_{k=1}^m(z-\zeta_{k,n})}
\prod_{j=1}^g(z-z_j).
$$
Теперь уже нетрудно проверить непосредственно, что эти
явные формулы сохраняют смысл, если при некотором
$j\in\{1,\dots,g\}$ точка~$z_j$ совпадает с~одним из краев
$j$-й лакуны~$L_j$. Каждую такую точку~$\mathbf z_j$ будем
считать как нулем функции $\psi^{(1)}(\mathbf x)$, так
и~нулем функции $\psi^{(2)}(\mathbf x)$,
$\mathbf x\in{\boldsymbol\Gamma}$. При этом соглашении
кусочно мероморфная на~$\mathfrak R$
функция~\eqref{eq67}--\eqref{eq68} дает решение
задачи~\ref{probR'} в~общем случае. Заметим, что так
как вес~$\rho$ -- голоморфная на ${\boldsymbol\Gamma}$
функция, то как функция $\psi(z^{(1)};n)$,
$z^{(1)}\in{D^{(1)}}$, так и функция $\psi(z^{(2)};n)$,
$z^{(2)}\in{D^{(2)}}$, голоморфно продолжаются
через~${\boldsymbol\Gamma}$ на другой лист римановой
поверхности. Поскольку, вообще говоря, $\rho\not\equiv1$,
то на~${\boldsymbol\Gamma}$ эти два голоморфных
продолжения не совпадают.

В~заключение отметим, что классическому случаю $g=0$ и
$\operatorname{supp}\mu=[-1,1]\cup\{\zeta_1,\dots,\zeta_m\}$
соответствуют
$$
\Pi_n(\mathbf z)=\prod_{k=1}^m\lambda(\zeta_k,\zeta_{k,n})\,
\frac{\varphi(z)-\varphi(\zeta_{k,n})}{z-\zeta_{k,n}}\,
\frac{z-\zeta_k}{1-\varphi(\zeta_k)\varphi(z)},
$$
где $\varphi(z)=z+\sqrt{z^2-1}$,
$\lambda(\zeta_k,\zeta_{k,n})\to1$, $n\to\infty$, и $D(z;\rho)$~--
функция Сегё (см.\ \cite{20} (\S\,5)).


\goodbreak

\end{document}